\documentclass[final, onefignum,onetabnum]{siamart190516}

\usepackage{lipsum}
\usepackage{amsfonts}
\usepackage{graphicx}
\usepackage{epstopdf}
\usepackage{algorithmic}
\usepackage{verbatim,amscd}
\usepackage{latexsym, bm, float, color}
\usepackage{mathtools}
\usepackage{tikz}
\usepackage{enumitem}

\ifpdf
  \DeclareGraphicsExtensions{.eps,.pdf,.png,.jpg}
\else
  \DeclareGraphicsExtensions{.eps}
\fi

\newsiamremark{remark}{Remark}
\newsiamremark{hypothesis}{Hypothesis}
\crefname{hypothesis}{Hypothesis}{Hypotheses}
\newsiamthm{claim}{Claim}

\newsiamthm{assumption}{Assumption}

\headers{Caputo Fractional Gradient Descent}{Y. Shin, J. Darbon, G. E. Karniadakis}

\title{A Caputo Fractional derivative-Based Algorithm for Optimization
	\thanks{Submitted to the editors DATE.
\funding{This work was funded by the PhILMS grant DE-SC0019453, the ARO MURI W911NF-15-1-0562 and the AFOSR MURI FA9550-20-1-0358.}}}

\author{Yeonjong Shin\thanks{Division of Applied Mathematics, Brown University, Providence, RI 02912
  (\email{yeonjong\_shin@brown.edu}, \email{jerome\_darbon@brown.edu}, \email{george\_karniadakis@brown.edu}).}
\and J\'er\^ome Darbon\footnotemark[2]
\and George Em Karniadakis\footnotemark[2]
\thanks{School of Engineering,
Brown University,
Providence, RI 02912, USA}
}

\usepackage{amsopn}

\newcommand{\R}{\mathbb{R}}
\newcommand{\x}{\textbf{x}}
\newcommand{\y}{\textbf{y}}
\newcommand{\w}{\textbf{w}}
\newcommand{\bc}{\textbf{c}}

\ifpdf
\hypersetup{
	pdftitle={Caputo Fractional Gradient Descent},
	pdfauthor={Y. Shin, J. Darbon, G. E. Karniadakis}
}
\fi

\begin{document}
	
	\maketitle
	
	\begin{abstract}
	    We propose a novel Caputo fractional derivative-based optimization algorithm.
        Upon defining the Caputo fractional gradient with respect to the Cartesian coordinate,
        we present a generic Caputo fractional gradient descent (CFGD) method.
		We prove that the CFGD yields the steepest descent direction 
		of a locally smoothed objective function.
		The generic CFGD requires three parameters to be specified,
		and
		a choice of the parameters yields a version of CFGD.
		We propose three versions -- non-adaptive, adaptive terminal and adaptive order.
		By focusing on quadratic objective functions, we provide a convergence analysis.
		We prove that the non-adaptive CFGD converges to 
		a Tikhonov regularized solution.
		For the two adaptive versions, 
		we derive error bounds, which show convergence to integer-order stationary point under some conditions.
		We derive an explicit formula of CFGD for quadratic functions.
        We computationally found that the adaptive terminal (AT) CFGD 
        mitigates the dependence on the condition number 
        in the rate of convergence and 
        results in significant acceleration over gradient descent (GD).
        For non-quadratic functions, 
        we develop an efficient implementation of CFGD using the Gauss-Jacobi quadrature,
        whose computational cost is approximately proportional to the number of the quadrature points and the cost of GD.
        Our numerical examples show that AT-CFGD results in acceleration over GD,
        even when a small number of the Gauss-Jacobi quadrature points (including a single point) is used.
	\end{abstract}
	
	
	
	\section{Introduction}

The gradient descent (GD) method to optimize a function 
dates back to Cauchy in 1800s \cite{Lemarechal_12_CauchyGD} and 
is one of the the most fundamental approaches in optimization.
It is an iterative algorithm to find a stationary point 
to an objective function. 
Due to its simplicity and scalability, 
GD and its variants have been widely used in many research fields,
in particular, machine learning \cite{Lecun_Nature15_DeepLearning, Ruder_16_GDoverview}.
Countless works have been devoted to GD-based methods, and the amount of literature is huge due to its importance.
Interested readers can consult the large number of textbooks
on the basics of optimization, e.g. \cite{bonnans.06.optimization,nesterov2004lectures,nocedalWright.06.book}.

Fractional calculus has been successfully employed 
in describing 
physical phenomena, e.g.,
anomalous transport, 
which classical models cannot capture (see \cite{d2020numerical,DElia_20_Unified} and references therein).
More recently, fractional calculus has been adapted in optimization algorithms \cite{Wei_20_FGD,Wang2017_NN17_FGD,Sheng2020convolutional,pu2013fractional,chen2017study}. 
Since fractional derivatives are extensions of integer-order derivatives,
one may naturally consider fractional gradient descent (FGD) as a generalization of gradient descent method. 
Many variants of FGD have been proposed
and shown to be effective in some applications \cite{Sheng2020convolutional,cheng2017innovative,khan2018fractional}.
Yet, many theoretical questions remain elusive.
For example, GD seeks to find an optimum by 
taking discrete steps in the direction of steepest descent.
However, which direction FGD follows
is not well understood.

We briefly review the existing literature on fractional calculus-based optimization methods. 
It has been pointed out in \cite{pu2013fractional,Wang2017_NN17_FGD}
that the set of stationary points of fractional gradient is different from the one of integer-order gradient.
In practice, integer-order stationary points are often sought.
To remedy this issue, \cite{Wei_20_FGD} proposed 
several heuristic variants of FGD.
These heuristics are useful for designing FGD algorithms,
yet no theoretical guarantees were provided.
%
While \cite{Wang2017_NN17_FGD} provided a convergence analysis, the results rely on multiple crude assumptions, such as uniform boundedness, which are not generally satisfied for many applications.
Another approach is to generalize the gradient flow, a continuum version of GD, to fractional time scale gradient flow \cite{Liang2020fractional,Hai2020gradient}.
This approach requires one to discretize the continuous flow appropriately to yield a numerical method.
Finally, we note that \cite{Nagaraj_20_Optimization} presents an abstract framework for convergence analysis of certain non-local calculus \cite{mengesha2015localization} based optimization. 

In this work, we propose a novel Caputo fractional gradient-based optimization algorithm, namely, the Caputo Fractional Gradient Descent (CFGD).
The Caputo fractional derivative \cite{Caputo1967linear} is one of the most popular fractional derivatives and is widely employed in modeling various physical phenomena, especially for initial value problems. 
Upon defining the Caputo fractional gradient with respect to the Cartesian coordinate \cite{tarasov2008fractional},
we define a generic CFGD algorithm \eqref{def:CFGD}.
The generic CFGD requires one to choose 
three parameters, and a choice of the parameters yields a version of CFGD.
The main findings are summarized as follows:
\begin{itemize}
    \item We prove that each direction generated by the generic CFGD is the steepest descent direction of a local smoothing of the original objective function (Theorem~\ref{thm:smoothing}).
    This answers the question of which direction CFGD follows towards minimizing the objective function.
    \item We propose three versions of CFGD -- non-adaptive, adaptive terminal and adaptive order.
    \item By focusing on quadratic objective functions, we provide a convergence analysis.
    \begin{enumerate}[label=(\roman*)]
        \item We prove that the non-adaptive CFGD converges linearly to
    a Tikhonov regularized solution,
    which is not a stationary point of integer-order gradient (Theorem~\ref{thm:fLSQ}).
        \item We provide error bounds of the two adaptive versions
    and show convergence to integer-order stationary point under some conditions
    (Theorems~\ref{thm:adapt-terminal} and~\ref{thm:adapt-order}).
    \end{enumerate}
    \item 
    For quadratic objective functions, 
    we derive an explicit formula of the Caputo fractional gradient (Corollary~\ref{cor:quad-CFGD}).
    We found that  
    the adaptive terminal CFGD effectively mitigates the dependence on the condition number in the convergence rate and yields a significant acceleration.
    \item For non-quadratic objective functions,
    based on Theorem~\ref{thm:formula-CFGD},
    we propose an efficient implementation of CFGD
    using the Gauss-Jacobi quadrature \eqref{eqn:Quad-CFGD}.
    We found that the adaptive terminal CFGD implemented by the proposed way 
    results in acceleration over GD
    even when a small number of the Gauss-Jacobi quadrature points (including a single point) is used.
\end{itemize}

\begin{figure}[!htbp]
	\centerline{
		\includegraphics[width=10.5cm]{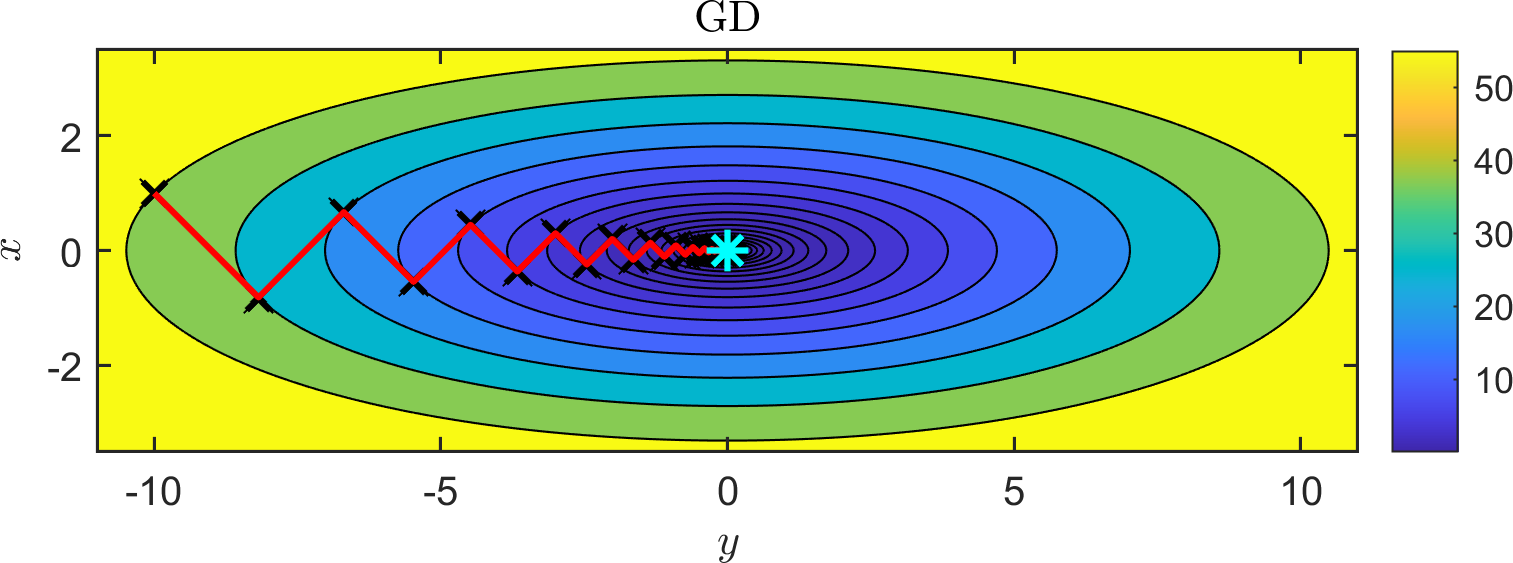}
	}
	\centerline{
		\includegraphics[width=10.5cm]{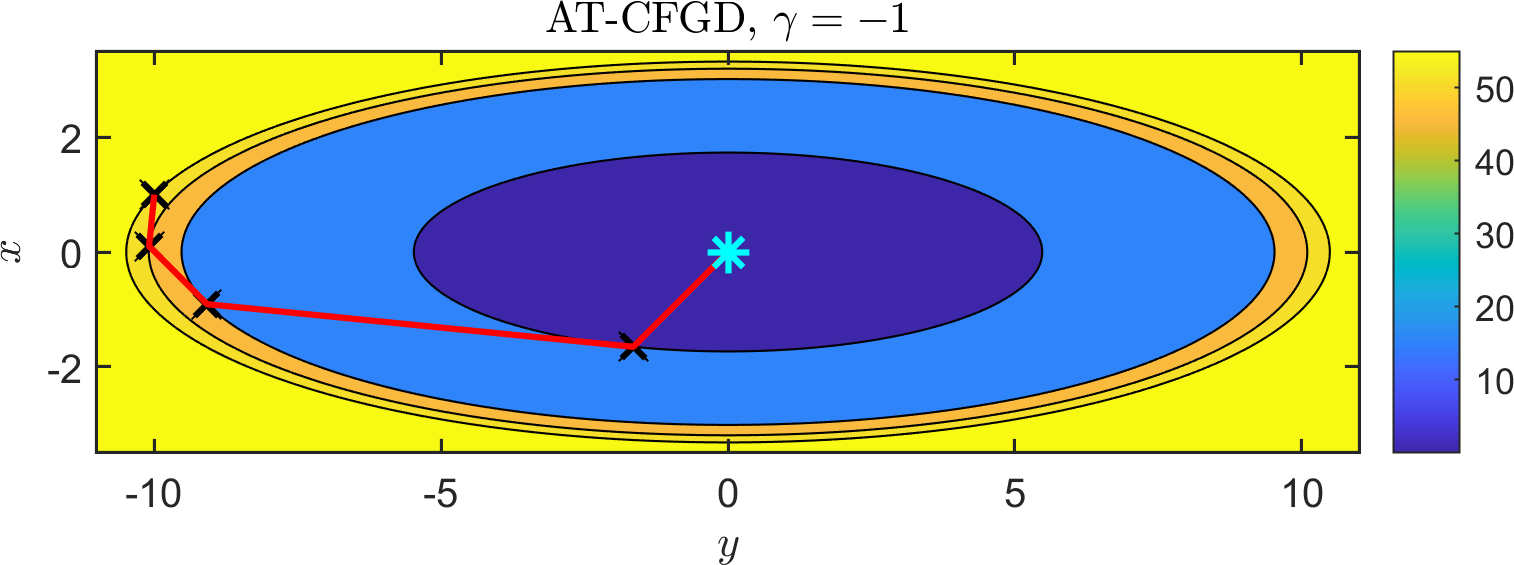}
	}
	\caption{A contour graph of the objective function $f(x,y)=10x^2+y^2$ (black) along with trajectories of (top) GD and (bottom) adaptive terminal CFGD.
	Each trajectory is shown as red solid line with the cross marks ($\times$).
	Both methods
	start at $\emph{\x}^{(0)} = (1,-10)$.
	}
	\label{fig:colormap}
\end{figure}

For quadratic objective functions, 
it is well known that GD converges linearly and the rate of convergence 
critically depends on
the condition number of the objective function.
A natural question is whether CFGD can  mitigate 
the dependence on the condition number
in the rate of convergence.
In Figure~\ref{fig:colormap}, we illustrate the performance of the adaptive terminal CFGD for a simple quadratic objective function $f(x,y) = 10x^2 + y^2$,
whose minimizer is the origin, marked as ($\ast$).
The trajectories of (top) GD and (bottom) CFGD are reported and they both start at $(1,-10)$.
We see that CFGD finds the optimal solution within the error of machine precision in merely four iterations.
GD behaves as expected showing 
a linear convergence.
However, since its rate depends on the condition number, more iterations are needed 
to reach the machine precision error.
Combined with more examples in Section~\ref{sec:example},
we found that  
the adaptive terminal CFGD effectively mitigates the dependence on the condition number in the convergence rate and yields a significant acceleration.

The rest of this paper is organized as follows. 
After describing the problem setup and introducing some preliminaries, 
the new Caputo fractional derivative-based algorithm is proposed in Section~\ref{sec:setup}.
A convergence analysis is presented  Section~\ref{sec:analysis}.
Numerical examples are provided in Section~\ref{sec:example} to verify our theoretical findings and to demonstrate 
the effectiveness of the proposed method.
	\section{Problem Setup and Method} \label{sec:setup}
We consider the general unconstrained minimization problem:
\begin{equation*}
	\min_{\x \in \mathbb{R}^d} f(\x),
\end{equation*}
where $f(\x)$ is a real-valued function.

The standard gradient descent method (GD) commences with an initial starting point $\x^{(0)}$
and updates the $k$-th iterated solution according to the following rule
\begin{equation*}
\x^{(k+1)} = \x^{(k)} - \eta_k\cdot \nabla_{\x} f(\x^{(k)}), \qquad k=0,1,\cdots
\end{equation*}
where $\eta_k > 0$ is the learning rate (or stepsize) at the $k$-th iteration. 
The method is theoretically and practically well-understood.
It is well-known that 
the GD converges linearly 
to a stationary point for many 
convex- and non-convex objective functions
as long as the learning rates are appropriately chosen \cite{bonnans.06.optimization,nesterov2004lectures,nocedalWright.06.book}.

\subsection{Caputo Fractional Derivative}
Since the fractional derivatives are not defined in a unified manner,
there exist multiple definitions.
In this paper, we focus on the fractional derivative in the sense of Caputo \cite{Caputo1967linear}
whose definition is given below.
\begin{definition}
	For $n \in \mathbb{N}\cup \{0\}$, let $\alpha \in (n-1,n]$.
	For $c, b \in \mathbb{R}$,
	let ${}_{c}I=[c,\infty)$ and $I_b = (-\infty, b)$
	be intervals.
	Let ${}_{c}\mathcal{D}^{\alpha}$ 
	and $\mathcal{D}^{\alpha}_{b}$ 
	be the sets of functions in $C({}_{c}I)$
	and $C(I_b)$, respectively,
	such that 
	\begin{equation*}
		\begin{split}
			{}_{c}\mathcal{D}^{\alpha} &= \{f \in C({}_{c}I): ~{}_{c}^{C}\!D_x^{\alpha} f \text{ exists and is finite in } {}_{c}I \}, \\
			\mathcal{D}^{\alpha}_b &= \{f \in C(I_b): ~{}_{x}^{C}\!D_{b}^{\alpha} f \text{ exists and is finite in } I_b \}.
		\end{split}
	\end{equation*}
	Here 
	${}_{c}^{C}\!D_x^{\alpha} f$ 
	is the left Caputo fractional derivative of $f$ of order $\alpha$ with the lower integral terminal $c$
	and
	$~{}_{x}^{C}\!D_{b}^{\alpha} f$ 
	is the right Caputo fractional derivative of $f$ of order $\alpha$ with the upper integral terminal $b$, defined respectively by
	\begin{equation*}
		\begin{split}
			(\text{Left}) 
			\quad 
			~{}_{c}^{C}\!D_x^{\alpha} f &:
			{}_{c}I \ni x \mapsto 
			\frac{1}{\Gamma(n-\alpha)}\int_c^x \frac{f^{(n)}(t)}{(x-t)^{\alpha-n +1}} dt \in \mathbb{R}, \\
			(\text{Right}) \quad ~{}_{x}^{C}\!D_b^{\alpha} f &:
			I_b \ni x \mapsto \frac{(-1)^n}{\Gamma(n-\alpha)}\int_x^b \frac{f^{(n)}(t)}{(t-x)^{\alpha-n +1}} dt \in \mathbb{R},
		\end{split}
	\end{equation*}
	where $f^{(n)}$ is the $n$-th derivative of $f$
	and $\Gamma$ is the Gamma function.
	For $c \in \mathbb{R}$, 
	we define a class $\mathcal{D}^{\alpha}(c)$ of functions which admit both left and right Caputo fractional derivatives with the integral terminal $c$:
	\begin{equation} \label{def:function-class}
	    \mathcal{D}^{\alpha}(c) = 
	    {}_{c}\mathcal{D}^{\alpha} \cap 
	    \mathcal{D}^{\alpha}_c.
	\end{equation}
\end{definition}

\begin{remark}
	For notational convenience, 
	${}_{b}^{C}\!D_x^{\alpha} f$
	is understood as 
	the right Caputo fractional derivative
	if $x < b$.
\end{remark}

\subsection{Smoothing Effect of Caputo Fractional Derivative}
Although the standard Caputo gradient descent method has been proposed and studied in several works \cite{Wei_20_FGD,Wang2017_NN17_FGD},
the motivation of using fractional derivatives
has been elusive in the context of optimization, 
except they are natural extensions of the integer order derivative.
In this section, we provide a mathematical justification of using the Caputo fractional derivative for optimization.

In the following theorem, 
we show that the Caputo fractional derivative 
induces an implicit regularization effect 
in the sense that it follows the steepest descent direction of a smoothing of the objective function $f$. 

\begin{theorem} \label{thm:smoothing}
    Let $f$ be a real-valued function defined on $\mathbb{R}$
    that admits a Taylor expansion around $c \in \mathbb{R}$.
    For $\alpha \in (0,1)$, $\beta, c \in \mathbb{R}$, 
    let ${}_{c}F_{\alpha,\beta}$ be a smoothing of $f$
    defined by
    \begin{equation*}
        {}_{c}F_{\alpha,\beta}(z) = f(c) + f'(c)(z-c) + \sum_{k=2}^{\infty} 
        C_{k,\alpha,\beta} \frac{f^{(k)}(c)}{k!}(z-c)^{k},
    \end{equation*}
    where $C_{k,\alpha,\beta} = \left(\frac{\Gamma(2-\alpha)\Gamma(k)}{\Gamma(k+1-\alpha)} + \beta\frac{\Gamma(2-\alpha)\Gamma(k)}{\Gamma(k-\alpha)}
        \right)$.
    Then, for any $x \ne c$, 
    the steepest descent direction of ${}_{c}F_{\alpha,\beta}$ at $x$ is 
    \begin{equation} \label{def:scaled-CFD}
       ({}_{c}F_{\alpha,\beta})'(x) = \frac{{}_{c}^{C}\!D^{\alpha}_x f}{{}_{c}^{C}\!D^{\alpha}_x I}(x)
       +\beta |x-c| \cdot \frac{{}_{c}^{C}\!D^{1+\alpha}_x f}{{}_{c}^{C}\!D^{\alpha}_x I}(x),
    \end{equation}
    where $I:\mathbb{R} \to \mathbb{R}$ is the identity map defined by $I(x) = x$.
\end{theorem}
\begin{proof}
    The proof can be found in Appendix~\ref{app:thm:smoothing}
\end{proof}

Theorem~\ref{thm:smoothing} indicates that 
the CFGD can be interpreted as
the steepest descent method applied on 
a particular smoothing ${}_{c}F_{\alpha,\beta}$ of $f$.
The parameters $\alpha, \beta$ control
the degree of smoothing through the coefficients $C_{k,\alpha,\beta}$.
If $\alpha =1$ and $\beta =0$, we recover the standard steepest descent direction of $f$.

\begin{figure}[!htbp]
	\centerline{
		\includegraphics[width=6.5cm]{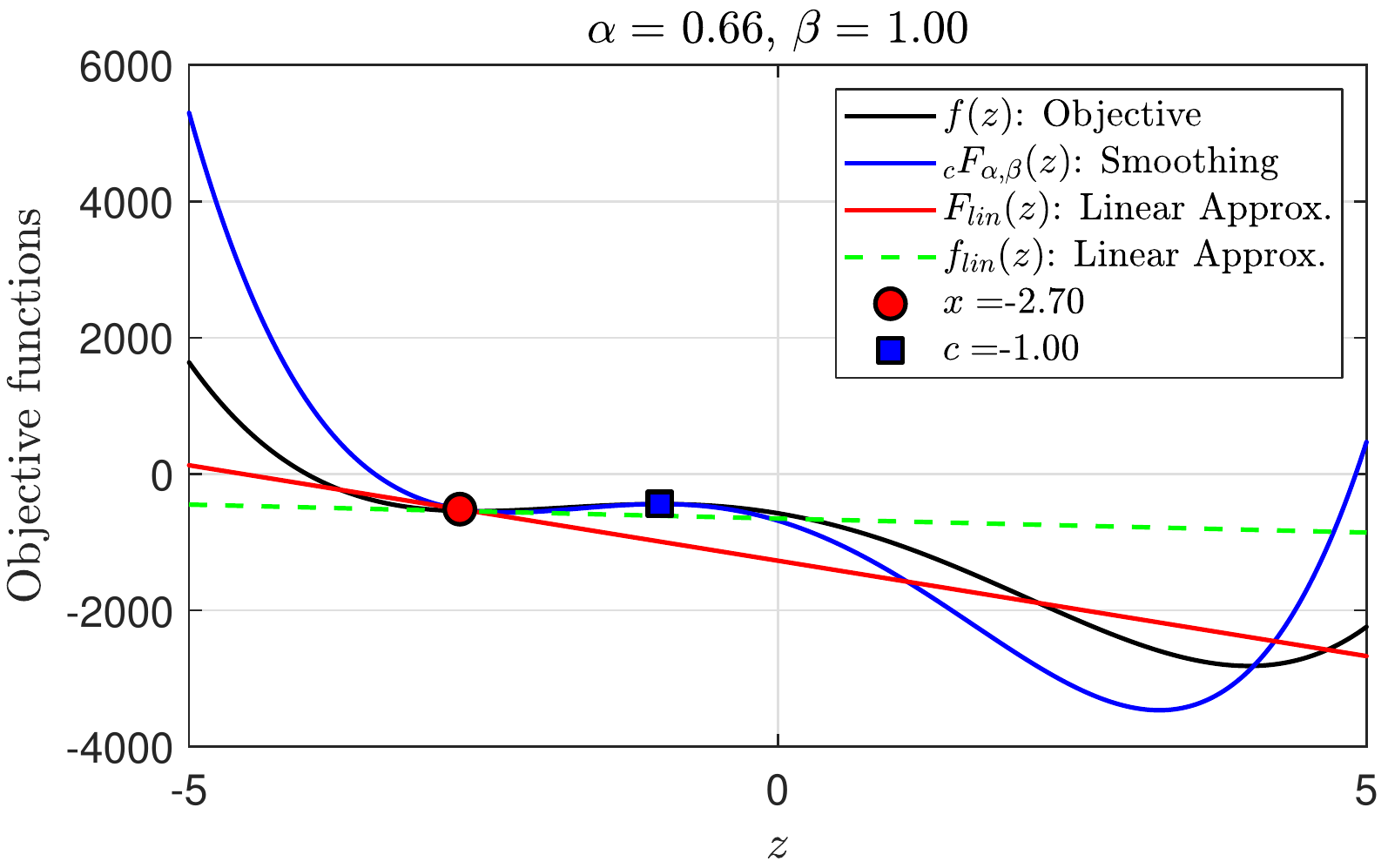}
		\includegraphics[width=6.5cm]{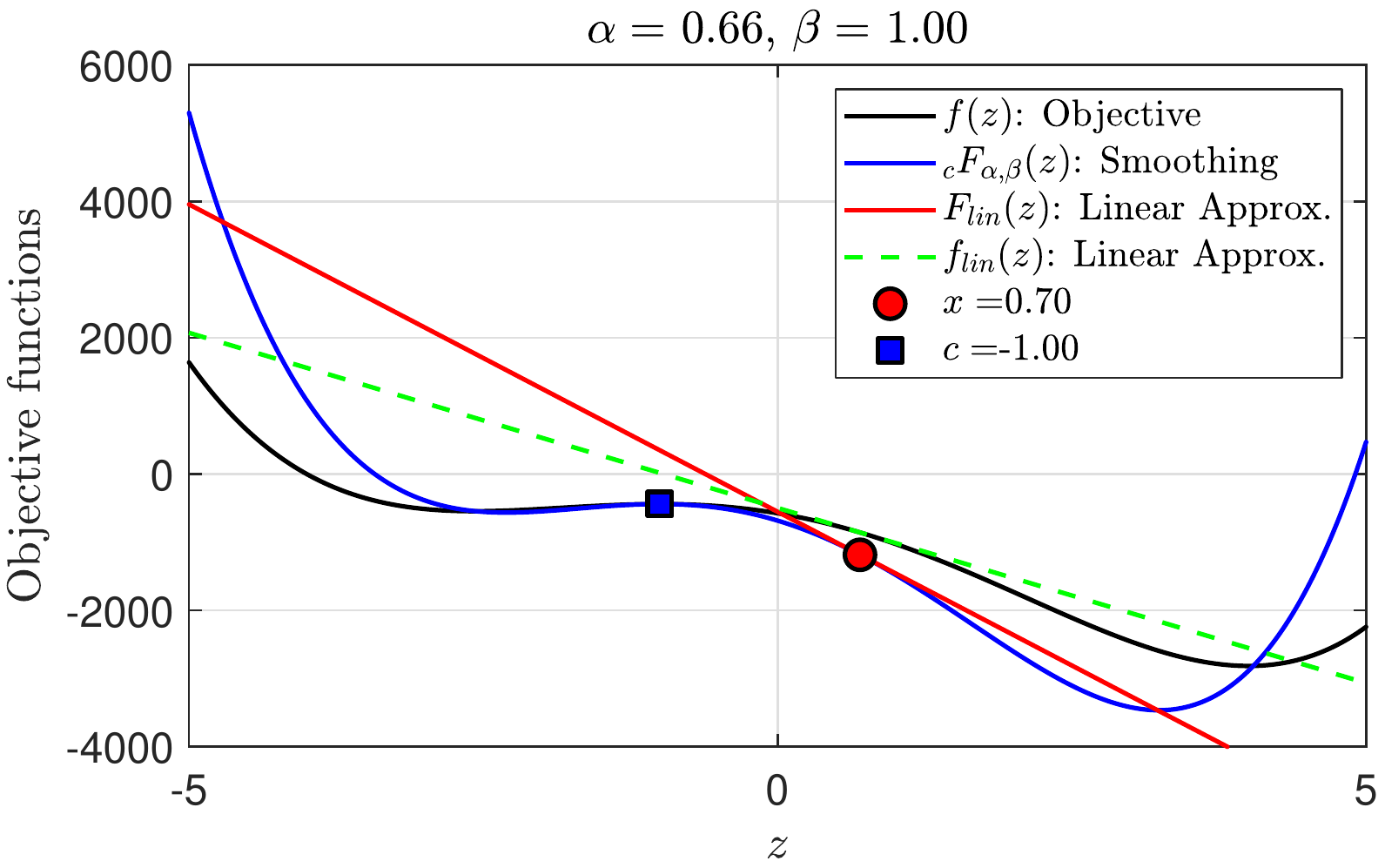}
	}
	\caption{The graphs of the objective function $f$ (black),
	its corresponding smoothing ${}_{c}F_{\alpha,\beta}$ (blue),
	the linear approximation $F_{\text{lin}}$ of ${}_{c}F_{\alpha,\beta}$ at $x$ (red)
	and the linear approximation $f_{\text{lin}}$ of $f$ at $x$ (green).
	Here $\alpha = 0.66$, $\beta = 1$
	and $c = -1$.
	(Left) $x = -2.7$ and 
	(Right) $x = 0.7$.
	}
	\label{fig:landscape}
\end{figure}

In Figure~\ref{fig:landscape},
we provide an illustration of the smoothing ${}_{c}F_{\alpha,\beta}$
for a particular objective function $f:\mathbb{R} \to \mathbb{R}$ defined by
\begin{equation*}
    f(z) = (z-6)(z+4)(7z^2+10z+24).
\end{equation*}
We set $\alpha = 0.66$, $\beta = 1$ and $c = -1$.
We also plot the linear approximation $f_{\text{lin}}$ of $f$
and the linear approximation $F_{\text{lin}}$ of ${}_{c}F_{\alpha,\beta}$ at $x = -2.7$ (left) and $x=0.7$ (right).
We note that the slope of $f_{\text{lin}}$ corresponds to the gradient of $f$ at $x$.
Similarly, the slope of $F_{\text{lin}}$ 
corresponds to the scaled Caputo fractional derivative \eqref{def:scaled-CFD} of $f$ at $x$.

\subsection{Caputo Fractional Gradient Descent}
We define the Caputo fractional gradient
by naturally extending the Caputo fractional derivative with respect to the Cartesian coordinate \cite{tarasov2008fractional}.

Let $f(\x)$ be a real-valued sufficiently smooth function defined on $\mathbb{R}^d$
with $\x = (x_1, \cdots, x_d)^\top \in \mathbb{R}^d$.
For $j=1,\dots, d$, let us define the functions $f_{j, \x}:\mathbb{\R} \to \mathbb{R}$ by 
$f_{j,\x}(y) = f(\x+ (y-x_j) e_j)$ where $e_j$ denotes the vector in $\mathbb{R}^d$ with a 1 in the $j$-th coordinate and 0's elsewhere.
For a vector $\textbf{c} = (c_1, \cdots,  c_d)^\top \in \mathbb{R}^d$,
following \cite{tarasov2008fractional}, 
we define the Caputo fractional gradient of $f$ by
\begin{equation} \label{def:Caputo-gradient}
    \prescript{C}{\bc}{\nabla_{\x}^{\alpha} f(\x)}:=
    \begin{bmatrix}
    {}_{c_1}^{C}\!D_{x}^{\alpha} f_{1,\x}(x_1)
    &\cdots
    &
    ~{}_{c_d}^{C}\!D_{x}^{\alpha} f_{d,\x}(x_d)
    \end{bmatrix}^\top \in \mathbb{R}^d,
\end{equation}
where $f_{j,\x}$ is assumed to be in $\mathcal{D}^{\alpha}(c_j)$ \eqref{def:function-class}
for all $j$.

We can now introduce a Caputo fractional gradient
descent method (CFGD) as follows:
Starting at an initial point $\x^{(0)}$,
the $k$-th iterated solution is updated by
\begin{equation*}
	\x^{(k+1)} = \x^{(k)} - \eta_k\cdot 
	~{}_{\textbf{c}_k}^{C}\!\nabla_{\x}^{\alpha} f(\x^{(k)}), 
	\qquad k=0,1,\dots
\end{equation*}
Here, $\eta_k$ is the learning rate of the $k$-th iteration
and $\textbf{c}_k$ is the lower/upper integral terminal at the $k$-th iteration.
We assume the range of $\alpha$ being $(0,1)$.


Motivated by Theorem~\ref{thm:smoothing},
we propose a properly scaled CFGD as follows.
Starting at an initial point $\x^{(0)}$,
the $k$-th iterated solution is updated by
\begin{equation} \label{def:CFGD}
\x^{(k+1)} = \x^{(k)} - \eta_k\cdot \vec{\textbf{d}}_k,
\qquad k=0,1,\dots,
\end{equation}
where $\x^{(k)} = (x^k_j)$, $\bc_k = (c^k_j)$, $\alpha_k \in (0,1)$, $\beta_k \in \mathbb{R}$ and 
\begin{equation} \label{def:direction}
    \vec{\textbf{d}}_k =
    \text{diag}\big({}_{c_j^k}^{C}\!D^{\alpha_k}_{x} I(x^k_j)\big)^{-1} 
    \bigg[ 
    {}_{\bc_k}^{C}\!\nabla^{\alpha_k}_{\x} f(\x^{(k)})
    +\beta_k \cdot \text{diag}\big( |x^k_j - c^k_j|\big) ~{}_{\bc_k}^{C}\!\nabla^{1+\alpha_k}_{\x} f(\x^{(k)}) \bigg].
\end{equation}
For a vector $\textbf{v} = (v_j) \in \mathbb{R}^d$,
either $\text{diag}(\textbf{v})$
or 
$\text{diag}(v_j)$
represents 
the diagonal matrix of size $d\times d$ whose $(j,j)$ component is $v_j$.
For the rest of the paper, 
we study the CFGD defined by \eqref{def:CFGD} and refer to it as the generic CFGD.

In general, computing fractional gradients 
is more expensive than computing integer-order gradients.
However, 
Theorem~\ref{thm:formula-CFGD}
provides an efficient way of evaluating 
$\vec{\textbf{d}}_k$
through the Gauss-Jacobi quadrature. 
\begin{theorem} \label{thm:formula-CFGD}
    Let $f(\emph{\x})$ be
	a real-valued sufficiently smooth function defined on $\R^d$
	and $f_{j,\emph{\x}}(y)$
	be the function defined in \eqref{def:Caputo-gradient}.
	Let $\alpha_k \in (0,1)$, 
	$\emph{\x}^{(k)} = (x_j^k)$,
	$\emph{\textbf{c}}_k = (c_j^k)$
	and $\Delta_j^k = |x_j^k - c_j^k|/2$.
	Then, 
	for $j=1,\dots, d$, we have 
	\begin{equation} \label{eqn:CFGD}
	    \begin{split}
	        \big(\vec{\emph{\textbf{d}}}_k \big)_j
		&= C_{\alpha_k} 
		\int_{-1}^{1}  f_{j,\emph{\x}}'(\Delta_j^k(1+u) + c_j^k)
		(1-u)^{-\alpha_k} du
		\\
		&\qquad+ C_{\alpha_k}
		\beta_k |x_j^k -c_j^k|
		\int_{-1}^{1}  f_{j,\emph{\x}}''(\Delta_j^k(1+u) + c_j^k)
		(1-u)^{-\alpha_k} du,
	    \end{split}
	\end{equation}
	where $C_{\alpha} = (1-\alpha)2^{-(1-\alpha)}$.
\end{theorem}
\begin{proof}
    The proof can be found in Appendix~\ref{app:thm:formula-CFGD}.
\end{proof}

We observe that \eqref{eqn:CFGD}
involves integrals 
that can be accurately evaluated by
the Gauss-Jacobi quadrature.
Let $\{(u_l,w_l)\}_{l=1}^s$ be 
the Gauss-Jacobi quadrature rule of $s$ points.
Then, $\vec{\textbf{d}}_k$ \eqref{eqn:CFGD} is approximated by $\vec{\textbf{d}}_{k,[s]}^{\text{QUAD}}$ \eqref{eqn:Quad-CFGD}
whose $j$-th component is defined by
\begin{equation} \label{eqn:Quad-CFGD}
    \begin{split}
        (\vec{\textbf{d}}_{k,[s]}^{\text{QUAD}})_j
        &= C_{\alpha_k}\sum_{l=1}^s w_l f_{j,\x}'(\Delta_j^k(1+u_l)+c_j^k)
        \\
        &\qquad 
        + 
        C_{\alpha_k}\beta_k|x_j^k-c_j^k| \sum_{l=1}^s w_l f_{j,\x}''(\Delta_j^k(1+u_l)+c_j^k).
    \end{split}
\end{equation}
By replacing $\vec{\textbf{d}}_k$
to $\vec{\textbf{d}}_{k,[s]}^{\text{QUAD}}$
in \eqref{def:CFGD},
we can efficiently implement the CFGD.
This allows the computation of fractional gradients to be embarrassingly parallel.
We briefly describe the complexity of a single iteration of the CFGD.
Suppose $s$ quadrature points are used
and the evaluations of $f_{j,\x}'$ and $f_{j,\x}''$ take $\mathcal{O}(K)$ FLOPS each.
Then it can be checked that 
computing $\vec{\textbf{d}}_{k,[s]}^{\text{QUAD}}$ requires $\mathcal{O}(sKd)$ FLOPS.
We note that a single iteration of the vanilla gradient descent takes $\mathcal{O}(Kd)$ FLOPS.
Thus, the computational complexity of CFGD grows only linearly in both the dimension $d$ and the number of quadrature points $s$.

\textit{Example.} 
In machine learning, 
the objective function is often defined through 
a set of training data
and a parametric model such as neural networks.
Suppose we have $m$ training data 
$\{(z_i,y_i)\}_{i=1}^m$, where $z_i, y_i \in \mathbb{R}$,
and use a two-layer neural network
$N(z;\x)= \sum_{j=1}^n a_{3,j}\phi(a_{1,j}z+a_{2,j})$
to fit these data.
Here $\phi$ is a nonlinear activation function
(e.g., $\tanh$ and sigmoid),
and 
$\x = \{a_{1,j},a_{2,j},a_{3,j}\}_{j=1}^n \in \mathbb{R}^{3n}$ is the set of network parameters
with $(\x)_l = a_{s,j}$ if $l = (s-1)n + j$, $s=1,2,3$, $j=1,\cdots,n$.
The goal is to find the optimal network parameter $\x^*$
that minimizes the following objective (loss) function:
\begin{equation*}
    f(\x) = \frac{1}{2}\sum_{i=1}^m \left(N(z_i;\x) - y_i \right)^2.
\end{equation*}
We discuss the computational complexity of
the evaluation of $\vec{\textbf{d}}_k^{\text{QUAD}}$ \eqref{eqn:Quad-CFGD}.
It suffices to calculate the complexities of $f_{j,\x}'$ and $f_{j,\x}''$, which turns out to be 
at most $K=\mathcal{O}(mn)$ FLOPS.
The detailed complexity calculation is reported in Appendix~\ref{app:complexity}.
Hence, if $s$ quadrature points are used
for the evaluation of $\vec{\textbf{d}}_k$ \eqref{eqn:CFGD},
the computational cost for the evaluation of \eqref{eqn:Quad-CFGD} is $\mathcal{O}(smn)$ FLOPS.

The aforementioned discussion implies that CFGD is roughly $s$-times more expensive than GD.
How large the number $s$ of quadrature points should be is an important question
to be addressed.
In Section~\ref{subsec:NNs},
we investigate the sensitivity of CFGD with respect to $s$.
While $\vec{\textbf{d}}_k^{\text{QUAD}}$ \eqref{eqn:Quad-CFGD} may no longer be an accurate approximation to $\vec{\textbf{d}}_k$ \eqref{eqn:CFGD}
for small $s$ values (including $s=1$), 
we empirically found that 
CFGD still results in acceleration over GD.
See Section~\ref{subsec:NNs} for more details.

\begin{remark}
    Another approach of evaluating fractional derivatives is to utilize a modern machine learning technique.
    In \cite{lu2019deeponet}, the authors demonstrated that 
    neural networks can learn linear and nonlinear operators. The resulting neural network is called a deep operator network (DeepONet). 
    We can pre-train a DeepONet for the purpose of evaluating Caputo fractional gradients 
    and utilize it in the CFGD algorithm.
    This deep learning approach will significantly lessen the computational cost of computing fractional gradients.
    We will pursue this direction in future study.
\end{remark}

\section{Convergence Analysis: Quadratic functions} \label{sec:analysis}

In this section, we present a convergence analysis of CFGD 
for the minimization of quadratic objective functions:
\begin{equation} \label{def:Quad-problem}
	\min_{\x} f(\x) = \frac{1}{2}\x^\top A\x + b^\top \x + c,
\end{equation}
where $\x, b \in \mathbb{R}^d, A = (a_{ij}) \in \mathbb{R}^{d\times d}$ and $c \in \mathbb{R}$.
Here, $A$ is assumed to be symmetric positive definite. 
It can be checked that the unique minimizer is
$\x^* = -A^{-1}b$.

As a special case, we consider the following regression problem.
Let 
\begin{equation*}
	W = \begin{bmatrix}
		\w_1, \w_2,\cdots, \w_m 
	\end{bmatrix} \in \mathbb{R}^{d \times m},
	\qquad
	\w_k = \begin{bmatrix}
		w_{k1} \\ \vdots \\ w_{kd}
	\end{bmatrix} \in \mathbb{R}^d,
	\quad 
	y = \begin{bmatrix}
		y_1 \\ \vdots \\ y_m
	\end{bmatrix} \in \mathbb{R}^m.
\end{equation*}
The least squares problem is formulated as follows
\begin{equation} \label{def:LSQ-problem}
	\min_{\x} f(\x) = \frac{1}{2} \|W^\top \x - y\|^2,
\end{equation}
where $\|\cdot\|$ is the Euclidean norm.
Assuming $W$ is of full rank with $m \ge d$, 
the least square solution is explicitly written as $\x^* = (WW^\top)^{-1}Wy$.
This is a special case of \eqref{def:Quad-problem}
with $A = WW^\top$ and $b = -Wy$.


For quadratic objective functions,
an explicit formula of the generic CFGD
can be obtained from Theorem~\ref{thm:formula-CFGD}.
\begin{corollary} \label{cor:quad-CFGD}
    Let the objective function $f(\emph{\x})$
    be a quadratic function of
    the form \eqref{def:Quad-problem}.
    Then, the direction defined in \eqref{def:direction} is explicitly given as 
    \begin{equation*}
        \vec{\emph{\textbf{d}}}_k = 
        A\emph{\x}^{(k)} + b 
	    + \gamma_{\alpha_k,\beta_k} \emph{\text{diag}}(\tilde{R})(\emph{\x}^{(k)}-\emph{\bc}_k),
    \end{equation*}
    where $\gamma_{\alpha,\beta} = \beta - \frac{1-\alpha}{2-\alpha}$
    and $\tilde{R} = (\sqrt{a_{11}},\cdots, \sqrt{a_{dd}})^\top$.
\end{corollary}
\begin{proof}
    The proof is readily followed from Theorem~\ref{thm:formula-CFGD}.
\end{proof}
Note that the integer-order gradient of 
the quadratic objective function at $\x^{(k)}$
is $A\x^{(k)} + b$.
Hence, the computational cost of CFGD
is roughly the same 
as the one of GD in this case.
Also, the fractional order $\alpha$ and the smoothing parameter $\beta$ depend only through $\gamma_{\alpha,\beta}$.
We often write $\gamma_{\alpha,\beta}$
as $\gamma$ if the context is clear.

Both GD and CFGD require one to 
determine appropriate learning rates (stepsizes).
For quadratic functions, 
an optimal learning rate could be obtained by 
the following line search:
\begin{align*}
    \min_{\eta} \frac{1}{2}(\x^{(k+1)})^\top A\x^{(k+1)} +b^\top \x^{(k+1)} + c,
\end{align*}
where $\x^{(k+1)} = \x^{(k)} - \eta_k \cdot \vec{\textbf{d}}_k$
for some direction vector $\vec{\textbf{d}}_k$ (either integer-order gradient or fractional gradient \eqref{def:direction}).
It then can be checked that the first-order optimality condition yields the optimal stepsize
\begin{equation} \label{def:LineSearch}
    \eta^*_k = \frac{\langle A\x^{(k)}+b, \vec{\textbf{d}}_k\rangle }{\vec{\textbf{d}}_k^\top A \vec{\textbf{d}}_k}.
\end{equation}

Depending on the choices of $\alpha_k, \beta_k, \bc_k$, 
the generic CFGD yields many variants.
In what follows, we propose and analyze three versions of CFGD.

\subsection{Non-adaptive Caputo Fractional Gradient Descent}
We consider the non-adaptive CFGD (NA-CFGD)
where we set 
$\alpha_k = \alpha$, $\beta_k = \beta$, and 
$\bc_k = \bc$ for all $k$
for some $\alpha \in (0,1)$, $\beta \in \mathbb{R}$ and $\bc \in \mathbb{R}^d$.

In Theorem~\ref{thm:fLSQ}, we will show that 
the stationary point of NA-CFGD 
is the solution to a Tikhonov regularization.
For reader's convenience and completeness,
we recall Tikhonov regularization \cite{golub1999tikhonov}:
\begin{equation*}
	\min_{\x} \|W^\top \x -y \|^2 + \gamma\|R^\top (\x-\overline{\x})\|^2,
\end{equation*}
where $R$ is some suitably chosen Tikhonov matrix.
It can be checked that 
the solution is given by
\begin{equation} \label{def:tik-opt}
	\x^*_{\text{Tik}}(\gamma) = 
	\overline{\x} + \left(WW^\top + \gamma RR^\top\right)^{-1}W(y-W^\top \overline{\x}),
\end{equation} 
assuming $WW^\top + \gamma RR^\top$ is invertible.

We are now in a position to 
present our convergence analysis of the non-adaptive CFGD.
\begin{theorem} \label{thm:fLSQ}
	Let the objective function $f$ have the form of 
	\eqref{def:LSQ-problem}.
	For $\alpha \in (0,1)$ and $\beta \in \mathbb{R}$, 
	let 
	$\tilde{A}_{\alpha,\beta}$ 
	be a matrix defined by
	\begin{equation*}
	    (\tilde{A}_{\alpha,\beta})_{ij} =
	\begin{cases}
	(\beta+\frac{1}{2-\alpha})\sum_{k=1}^m w_{ki}^2 & \text{if } i = j \\
	\sum_{k=1}^m w_{ki}w_{kj}	  & \text{if } i \ne j
	\end{cases}.
	\end{equation*}
	Suppose $\tilde{A}_{\alpha,\beta}$ 
	is a positive definite matrix
	and $\sigma_{\max}$ is its largest singular value.
	Let $\alpha_k = \alpha$, $\beta_k = \beta$, $\emph{\bc}_k = \emph{\bc}$ and 
	$\eta_k = \frac{\eta}{\sigma_{\max}}$ for some $0 < \eta < 2$ and $\emph{c} \in \mathbb{R}^d$.
	Then, the $k$-th iterated solution of \eqref{def:CFGD}
	satisfies 
	\begin{equation*}
		\|\emph{\x}^{(k)} - \emph{\x}^*_{\emph{\text{Tik}}}(\gamma)\|^2 \le \|\emph{\x}^{(0)} - \emph{\x}^*_{\emph{\text{Tik}}}(\gamma)\|^2 |1- \frac{\eta}{\kappa_{\alpha,\beta}}|^k,
	\end{equation*}
	where 
	$\kappa_{\alpha,\beta}$ is the condition number of $\tilde{A}_{\alpha,\beta}$
	and 
	$\emph{\x}^*_{\emph{\text{Tik}}}$ is the solution to the Tikhonov regularization \eqref{def:tik-opt}
	with $\overline{\emph{\x}} = \emph{\textbf{c}}$,
	$\gamma = \beta-\frac{1-\alpha}{2-\alpha}$
	and 
	$
	R = \emph{\text{diag}}(\sqrt{\sum_{k=1}^m w_{kj}^2}). 
	$
\end{theorem}
\begin{proof}
	The proof can be found in Appendix~\ref{app:thm:fLSQ}.
\end{proof}

Theorem~\ref{thm:fLSQ} shows that NA-CFGD converges linearly to the stationary point, which is the solution to a certain Tikhonov regularization.
The convergence rate depends only on $A$ and
$\beta + \frac{1}{2-\alpha}$.
We note that when $\alpha = 1$ and $\beta = 0$,
it can be checked from Corollary~\ref{cor:quad-CFGD}
that we recover gradient descent. 
Also, 
if $\beta + \frac{1}{2-\alpha} > 1$, 
NA-CFGD expects to converge faster than
GD as 
the condition number of $\kappa_{\alpha,\beta}$ is smaller than $\kappa_{1,0}$.
However, the fractional stationary point is not the same as the integer-order stationary point, which is often sought in practice.

With the goal of finding the integer-order stationary point,
in the following two subsections, 
we consider two adaptive versions of CFGD.
One is adaptive terminal 
and the other is adaptive order.

\subsection{Adaptive Terminal Caputo Fractional Gradient Descent}
We consider the adaptive terminal CFGD (AT-CFGD),
which uses adaptive integral terminal $\bc_k = \x^{(k-L)}$ for some positive integer $L$,
while 
$\alpha_k = \alpha$ and $\beta_k = \beta$ for all $k$
for some constants $\alpha \in (0,1)$, $\beta \in \mathbb{R}$.
Hence, AT-CFGD commences with $L$ initial points $\{\x^{(-j)}\}_{j=0}^L \subset \mathbb{R}^d$. 



The following theorem shows 
an error bound of AT-CFGD with respect to 
the optimum $\x^*$ of \eqref{def:Quad-problem}.

\begin{theorem} \label{thm:adapt-terminal}
	Let the objective function $f$ have the form of 
	\eqref{def:Quad-problem}
	and $\emph{\x}^*$ be the corresponding optimal solution.
	For $\alpha \in (0,1), \beta \in \mathbb{R}$ and $L \in \mathbb{N}_{\ge 1}$, 
	let $\gamma_{\alpha,\beta} = \beta - \frac{1-\alpha}{2-\alpha}$ and 
	$\{\emph{\x}^{(-j)}\}_{j=0}^L \subset \mathbb{R}^d$
	be initial points.
	Suppose 
	$\alpha_k = \alpha$, $\beta_k = \beta$,
	$\eta_k = \eta$,  
	$\emph{\bc}_k = \emph{\x}^{(k-L)}$.
	Then, the $k$-th iterated solution of \eqref{def:CFGD}
	satisfies 
	\begin{equation*}
		\|\emph{\x}^{(k)} - \emph{\x}^*\| \le 
		\sum_{j=0}^L \|\mathcal{A}_{k,j}\|\|\emph{\x}^{(-j)} - \emph{\x}^*\|,
	\end{equation*}
	where $\mathcal{A}_{k,j} \in \mathbb{R}^{d\times d}$ are matrices defined recursively by
	\begin{align*}
	    \mathcal{A}_{k,0} = \mathcal{A}_{k-1,0}\mathcal{A}_{1,0} + \mathcal{A}_{k-1,1}, \quad
	    \mathcal{A}_{k,L} = \mathcal{A}_{k-1,0}\mathcal{A}_{1,L},
	    \quad
	    \mathcal{A}_{k,j} = \mathcal{A}_{k-1,j+1},
	\end{align*}
	for $1 \le j < L$ and $k=2,\dots$ starting with 
	$\mathcal{A}_{1,0} = I - \eta (A + \gamma_{\alpha,\beta}\emph{\text{diag}}(A))$,
	$\mathcal{A}_{1,L} = \eta  \gamma_{\alpha,\beta}\emph{\text{diag}}(A)$
	and $\mathcal{A}_{1,j} = 0$ for $1 \le j < L$.
	The matrix norm $\|\mathcal{A}_{k,j}\|$
	is understood as the spectral norm.
	Furthermore,
	if $\lim_{k\to \infty} \|\mathcal{A}_{k,0}\|= 0$,
	we have
	\begin{equation*}
	    \lim_{k\to \infty} \|\emph{\x}^{(k)} - \emph{\x}^*\| = 0.
	\end{equation*}
\end{theorem}
\begin{proof}
    The proof can be found in Appendix~\ref{app:thm:adapt-terminal}.
\end{proof}

Although Theorem~\ref{thm:adapt-terminal} provides 
an error bound with respect to the integer-order stationary point, 
further investigation 
is needed in understanding the dynamics of $\mathcal{A}_{k,0}$ for convergence. 
Yet, we empirically found that AT-CFGD not only 
converges but also converges significantly faster than GD.
We remark that 
the convergence rate of AT-CFGD
does not explain 
such a significant acceleration.
As shown in Figure~\ref{fig:colormap},
the first two iterations of AT-CFGD 
do not show a significant improvement (actually no better than GD).
This indicates that linear convergence (uniform rate) 
does not explain the accelerated convergence of AT-CFGD.
More numerical tests are reported in Section~\ref{sec:example}, while we provide details of Figure~\ref{fig:colormap} below.

\begin{figure}[!htbp]
	\centerline{
		\includegraphics[width=6cm]{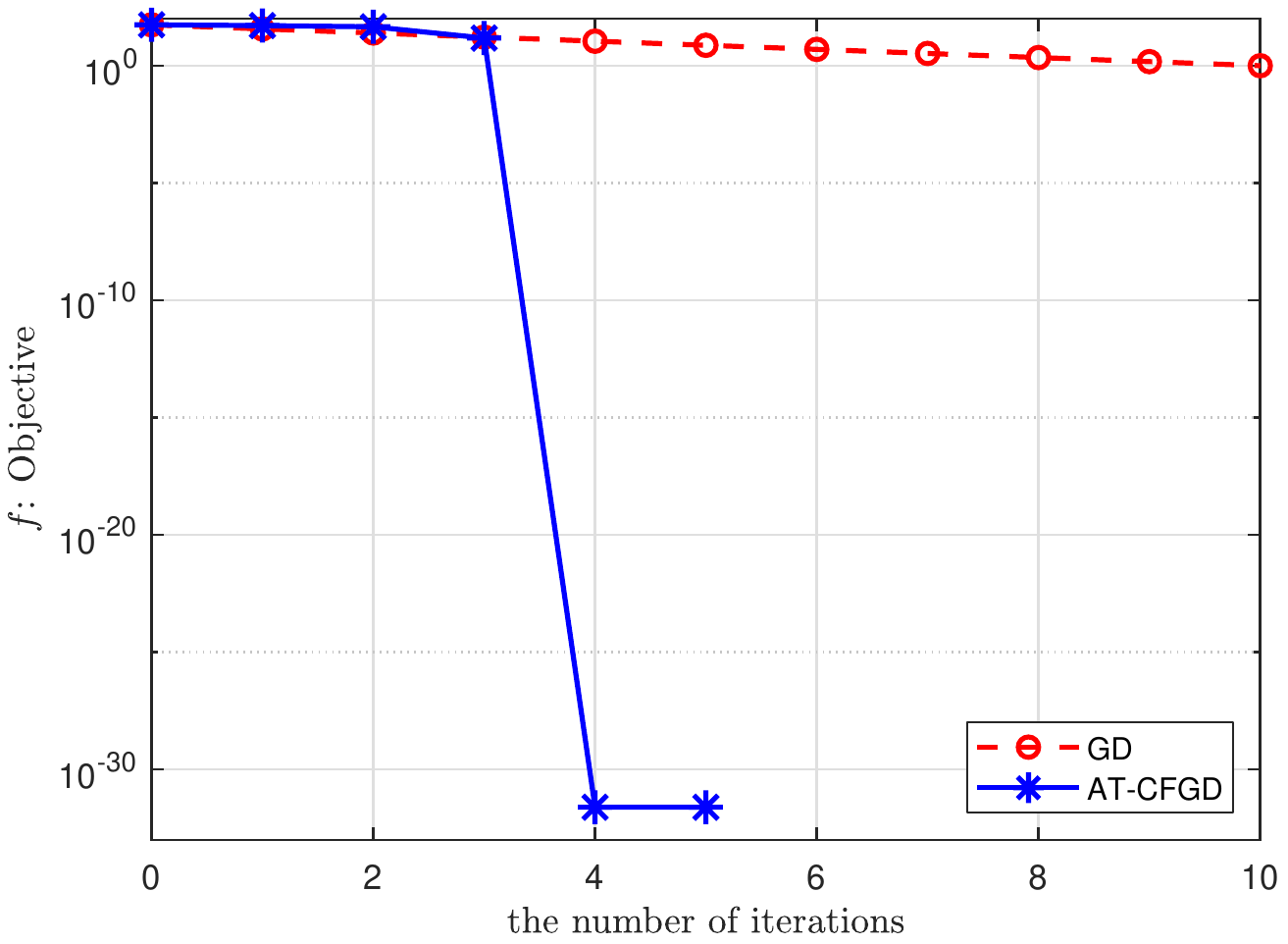}
	}
	\caption{
	The same example of Figure~\ref{fig:colormap} is considered.
	The objective function is $f(x,y) = 10x^2 + y^2$.
	The objective values versus the number of iterations.
	AT-CFGD finds the global minimum within merely 4 iterations, while GD requires a large number of iterations to reach the same level.
	}
	\label{fig:loss-convg}
\end{figure}

In Figure~\ref{fig:colormap}, we employ AT-CFGD with $L = 1$, $\gamma = -1$,
$\x^{(-1)} = (-1,-1)$, and the optimal stepsize of \eqref{def:LineSearch}.
Note that $\alpha$ and $\beta$ depend only through $\gamma_{\alpha,\beta}$ (Corollary~\ref{cor:quad-CFGD}).
We further report the objective values versus the number of iterations 
in Figure~\ref{fig:loss-convg}.
We see that (as expected) GD converges linearly to the optimum whose rate depends on the condition number of the model matrix.
While AT-CFGD does not exhibit a linear convergence, 
it finds the optimal solution within the error of machine accuracy in merely four iterations.

\subsection{Adaptive Order Caputo Fractional Gradient Descent}
We consider the adaptive order CFGD (AO-CFGD), which uses adaptive fractional order $\alpha_k$,
adaptive smoothing parameter $\beta_k$,
while $\bc_k = \bc$ for some $\bc \in \mathbb{R}^d$ for all $k$.

Let $\{(\alpha_s,\beta_s)\}_{s\ge 1}$ be a sequence such that 
$\alpha_s \in (0,1]$ and $\beta_s \in [\frac{1-\alpha_s}{2-\alpha_s},\infty)$.
Let $\{\gamma_s\}_{s\ge 1}$ be a nonnegative sequence defined by 
\begin{equation} \label{def:adapt-order}
    \gamma_s = \beta_s -\frac{1-\alpha_s}{2-\alpha_s}.
\end{equation}
Let $\{k_s\}_{s \ge 1}$ be a sequence of positive integers.
For a positive integer $s$, 
AO-CFGD consists of $s$ stages.
The first stage starts with an initial point $\x^{(0)}$
and apply the CFGD of \eqref{def:CFGD} with
$\alpha_1, \beta_1, \bc$
for the first $k_1$ iterations.
Let us denote the $k_1$-th iterated solution by $\x^{(k_1)}_1$.
For $s \ge 2$,
the $s$-th stage starts with $\x^{(k_{s-1})}_{s-1}$
and apply the CFGD with 
$\alpha_s, \beta_s, \bc$ for the next $k_s$ iterations.
The resulting solution is denoted by $\x^{(k_s)}_s$.

In Theorem~\ref{thm:adapt-order}, 
we provide an error bound of AO-CFGD
with respect to the optimal solution $\x^*$.

\begin{theorem} \label{thm:adapt-order}
	Let $\{\gamma_s\}$ be a nonnegative convergent sequence to 0 defined in
	\eqref{def:adapt-order}.
	Let $\kappa_{s}$ be the condition number of 
	$\tilde{A}_{\alpha_s,\beta_s}$ defined in Theorem~\ref{thm:fLSQ}.
	Let $\eta_{s,k}$ be the learning rate of the $k$-th iteration at the $s$-th stage.
	Let $\eta_{s,k} = \eta_s$
	for some $\eta_s > 0$.
	Let
	$R_s = |1 - \eta_s/\kappa_s|^{k_s/2}$ for all $s$.
	Then, the solution to AO-CFGD after $s$ stages satisfies 
	\begin{align*}
	    &\|\emph{\x}_s^{(k_s)} - \emph{\x}^*\|
	    \\
	    &\le 
	    \prod_{j=0}^{s-1} R_{s-j}
	    \|\emph{\x}^{(0)} - \emph{\x}^*_{\emph{\text{Tik}}}(\gamma_1)\|
	    +C\left\{ \sum_{i=1}^{s-1} \left(\prod_{j=0}^{i-1} R_{s-j}\right) |\gamma_{s-i} - \gamma_{s-i+1}| + |\gamma_s|
	    \right\},
	\end{align*}
	where $\emph{\x}^*_{\emph{\text{Tik}}}$ is the solution to the Tikhonov regularization \eqref{def:tik-opt}
	and $C$ is a constant defined in Appendix~\ref{app:thm:adapt-order}
	that depends only on $W$, $R$, $\emph{\x}^*$, $\emph{\bc}$
	and the range of $\{\gamma_s\}$.
	Furthermore, if $k_s = k$ and $R_s < \rho^k$ for some $\rho \in (0,1)$ for all $s$,
	we have
	\begin{align*}
	    \lim_{k \to \infty} \|\emph{\x}_s^{(k)} - \emph{\x}^*\|
	    \le C|\gamma_s|, \quad
	    \lim_{s \to \infty} \|\emph{\x}_s^{(k)} - \emph{\x}^*\|
	    \le C \frac{\rho^k}{1-\rho^k},
	    \quad
	    \lim_{s,k \to \infty} \|\emph{\x}_s^{(k)} - \emph{\x}^*\| = 0.
	\end{align*}
\end{theorem}
\begin{proof}
	The proof can be found in Appendix~\ref{app:thm:adapt-order}.
\end{proof}

We note that the idea of using 
adaptive fractional order
appeared in \cite{Wei_20_FGD},
where multiple heuristic adaptation strategies
were presented. 
While some promising empirical results were reported in \cite{Wei_20_FGD}, 
finding an optimal strategy requires more investigation
and remains a challenging problem.
Due to these reasons, 
we do not consider 
AO-CFGD in numerical tests in Section~\ref{sec:example}.

	\section{Numerical Examples}
\label{sec:example}

We present numerical examples to
verify our theoretical findings and demonstrate the performance of our proposed Caputo fractional gradient descent.

The generic CFGD \eqref{def:CFGD} involves a set of hyperparameters to be chosen -- 
the fractional order $\alpha_k$,
the smoothing parameter $\beta_k$
and the integral terminal $\bc_k$. 
We focus on two versions -- non-adaptive CFGD, which is referred to as NA-CFGD
and adaptive terminal CFGD,
which is referred to as AT-CFGD.
We recall that NA-CFGD sets all the parameters to be constants. 
AT-CFGD sets $\alpha_k = \alpha, \beta_k = \beta$ and $\bc_k = \x^{(k-L)}$
for some $\alpha \in (0,1), \beta \in \mathbb{R}$ and $L \in \mathbb{N}_{\ge 1}$. 
Unless otherwise stated, 
we employ the optimal learning rate (stepsize) of \eqref{def:LineSearch}
in both GD and CFGD.

\subsection{Quadratic Objective Function: Random data}
We consider quadratic objective functions
of the form of \eqref{def:LSQ-problem}.
We generate a matrix $W$ of size $d \times m$ by 
sampling each entry independently from a normal distribution $\mathcal{N}(0,1/m)$.
Similarly, we randomly generate the vector $y$ of size $m$
from $\mathcal{N}(0,I_m)$, where $I_m$ is the identity matrix of size $m$.
In terms of \eqref{def:Quad-problem}, 
the model matrix $A$ is $WW^\top$
and $b$ is $-Ay$.

\begin{figure}[!htbp]
	\centerline{
		\includegraphics[width=6.2cm]{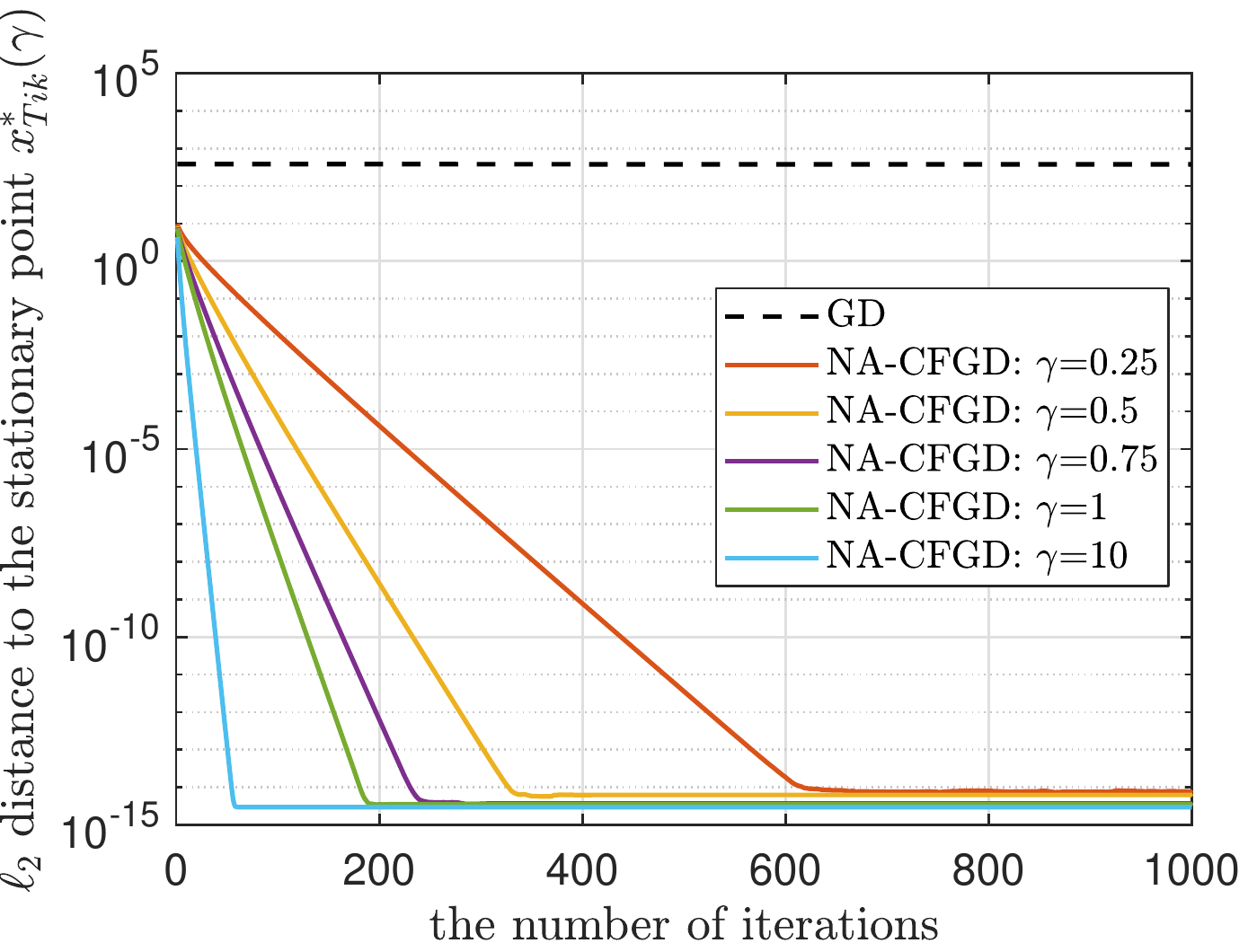}
		\includegraphics[width=6.2cm]{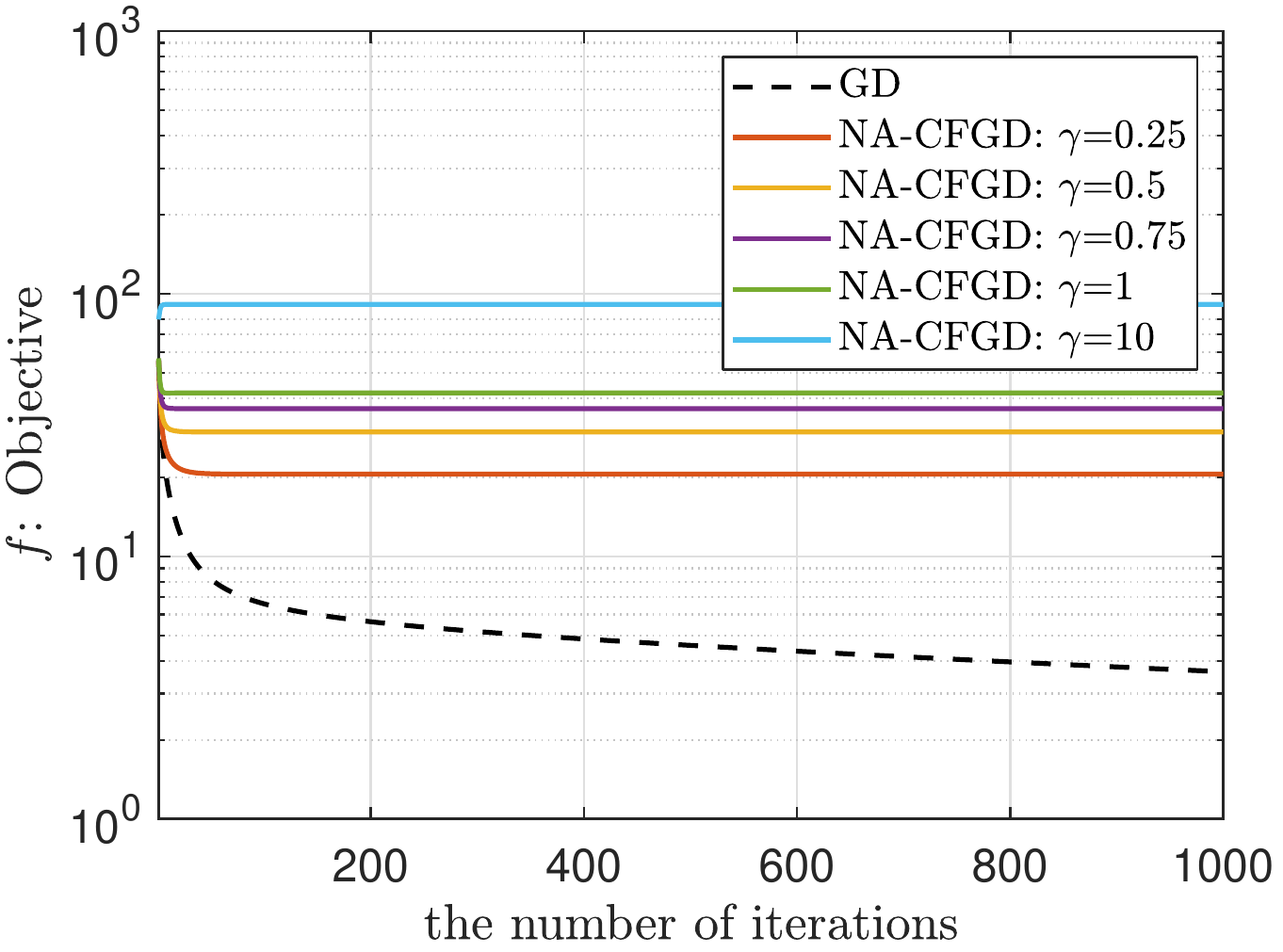}
	}
	\caption{(To be viewed in color)  
	(Left) The $\ell_2$ distance to the optimal point $\emph{\x}^*_{\emph{\text{Tik}}}$ of 
	Tikhonov regularization 
	with respect to the number of iterations at varying $\gamma_{\alpha,\beta}$.
	(Right) The objective values versus the number of iterations.
	SD-CFGD does not decrease the objective function
	as it does not converge to the stationary point $\emph{\x}^*$ of the objective. 
	}
	\label{fig:D100-Tikhonov}
\end{figure}

We first verify the convergence of NA-CFGD 
to the solution to Tikhonov regularization.
We set the integral terminal vector $\bc$ to a vector whose elements are all ones
and the initial starting point $\x^{(0)}$ to a random vector from the standard normal distribution.
In Figure~\ref{fig:D100-Tikhonov}, we show the results for $d=m=100$. 
On the left, the $\ell_2$ distance to the solution $\x^*_{\text{Tik}}(\gamma)$ \eqref{def:tik-opt} is reported 
with respect to the number of iterations at varying $\gamma \in \{0.15, 0.25, 0.5, 0.75, 1, 10\}$.
Note that gradient descent corresponds to the case of $\gamma =0$.
As expected from Theorem~\ref{thm:fLSQ}, 
we clearly see that NA-CFGD converges to $\x^*_{\text{Tik}}$
at a significantly faster rate compared to those of GD.
This is not a surprise as Tikhonov regularization induces 
a smaller condition number of the model matrix,
that results in a fast convergence.
At the same time, it shows that NA-CFGD does not find the integer-order stationary point.
On the right, the objective values are reported. 
As expected, NA-CFGD does not decrease the objective function,
while it already reaches to its own stationary point. 
We note that in this case, the model matrix $A = WW^\top$ is ill-conditioned, whose condition number is greater than $10^5$.
This explains an extremely slow convergence by GD.

Next, we demonstrate the performance of adaptive terminal CFGD, which enables to reach the integer-order stationary point $\x^*$.
We consider an ill-conditioned model matrix generated by the following procedure.
First, we generate a random matrix $W_0$ according to the aforementioned manner.
After computing a singular value decomposition $W_0 = US_0V^\top$,
we consider the matrix $S$ obtained from $S_0$ by replacing 
the largest singular value to $10$
and the smallest singular value to $0.1$.
We then obtain a matrix defined by $W = USV^\top$ whose condition number is at least $10^2$.
Therefore, the condition number of the model matrix $A = WW^\top$ is at least $10^4$.
In what follows, we set $d=m=20$.
For both GD and AT-CFGD, $\x^{(0)}$ is randomly chosen from the uniform distribution of $[-1,1]^d$.
The initial points $\{\x^{(-j)}\}_{j=1}^{L}$ for AT-CFGD
are randomly chosen from the standard normal distribution.

\begin{figure}[!htbp]
	\centerline{
		\includegraphics[width=7.5cm]{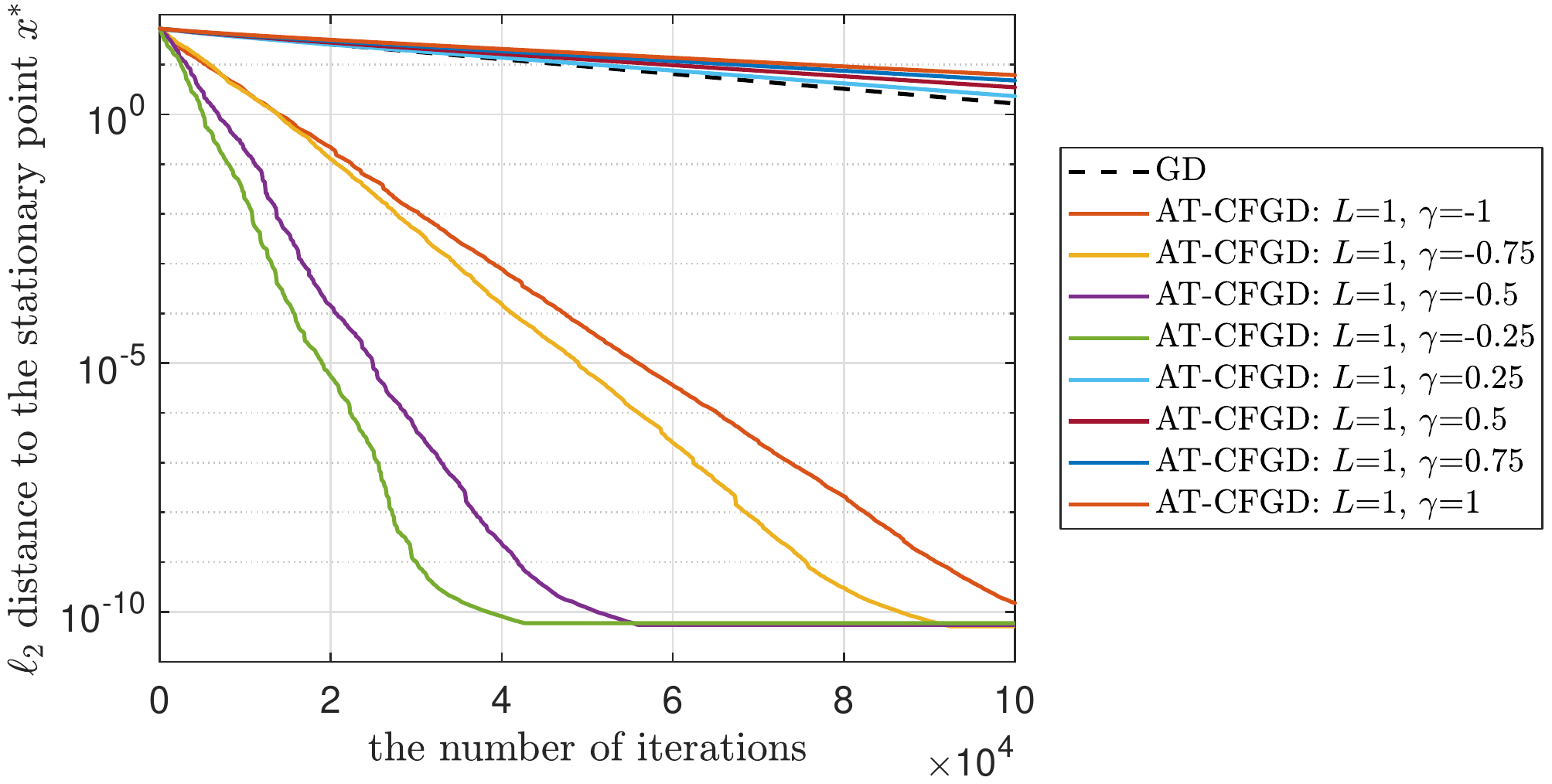}
	}
	\centerline{
		\includegraphics[width=4.3cm]{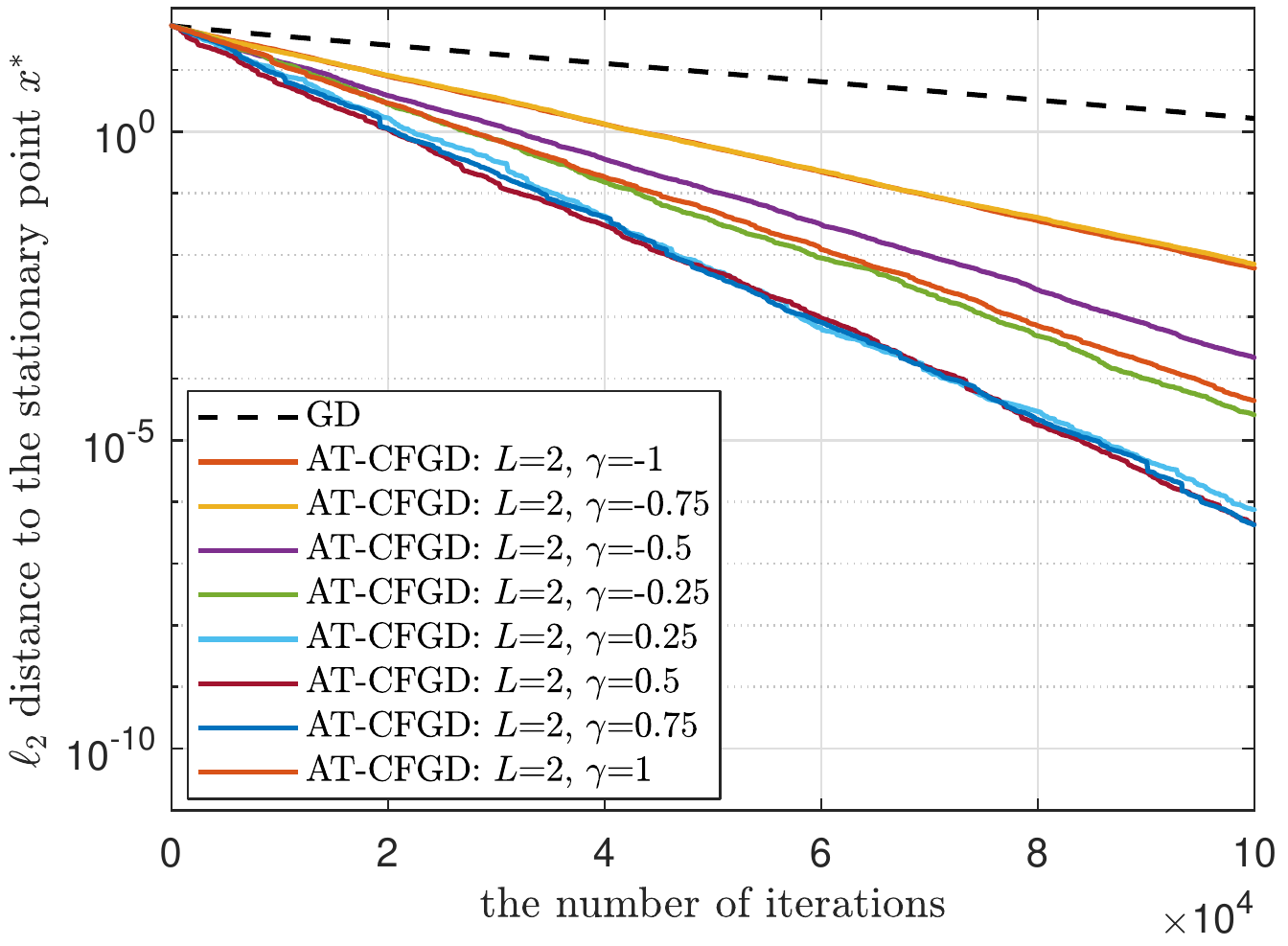}
		\includegraphics[width=4.3cm]{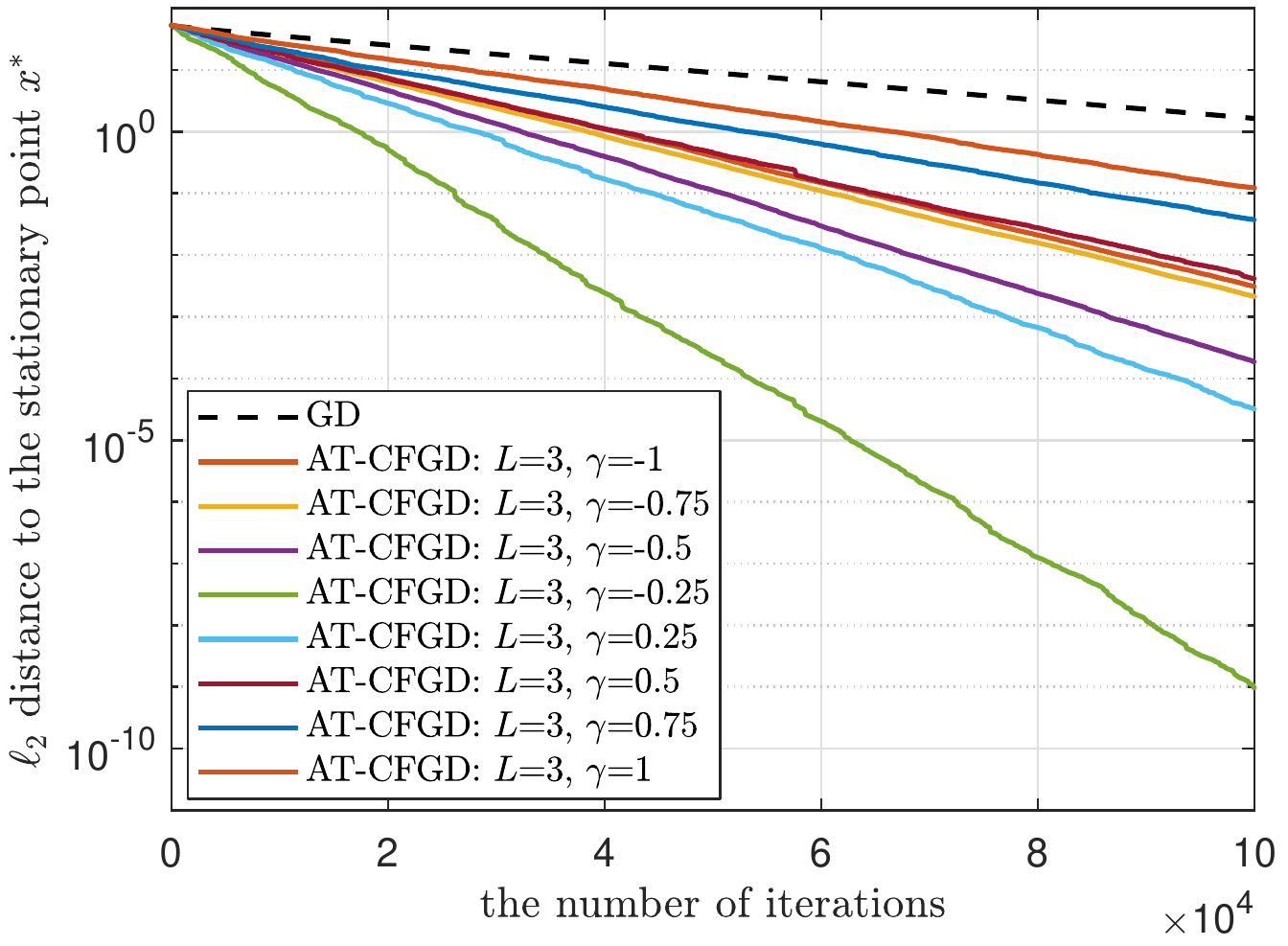}
		\includegraphics[width=4.3cm]{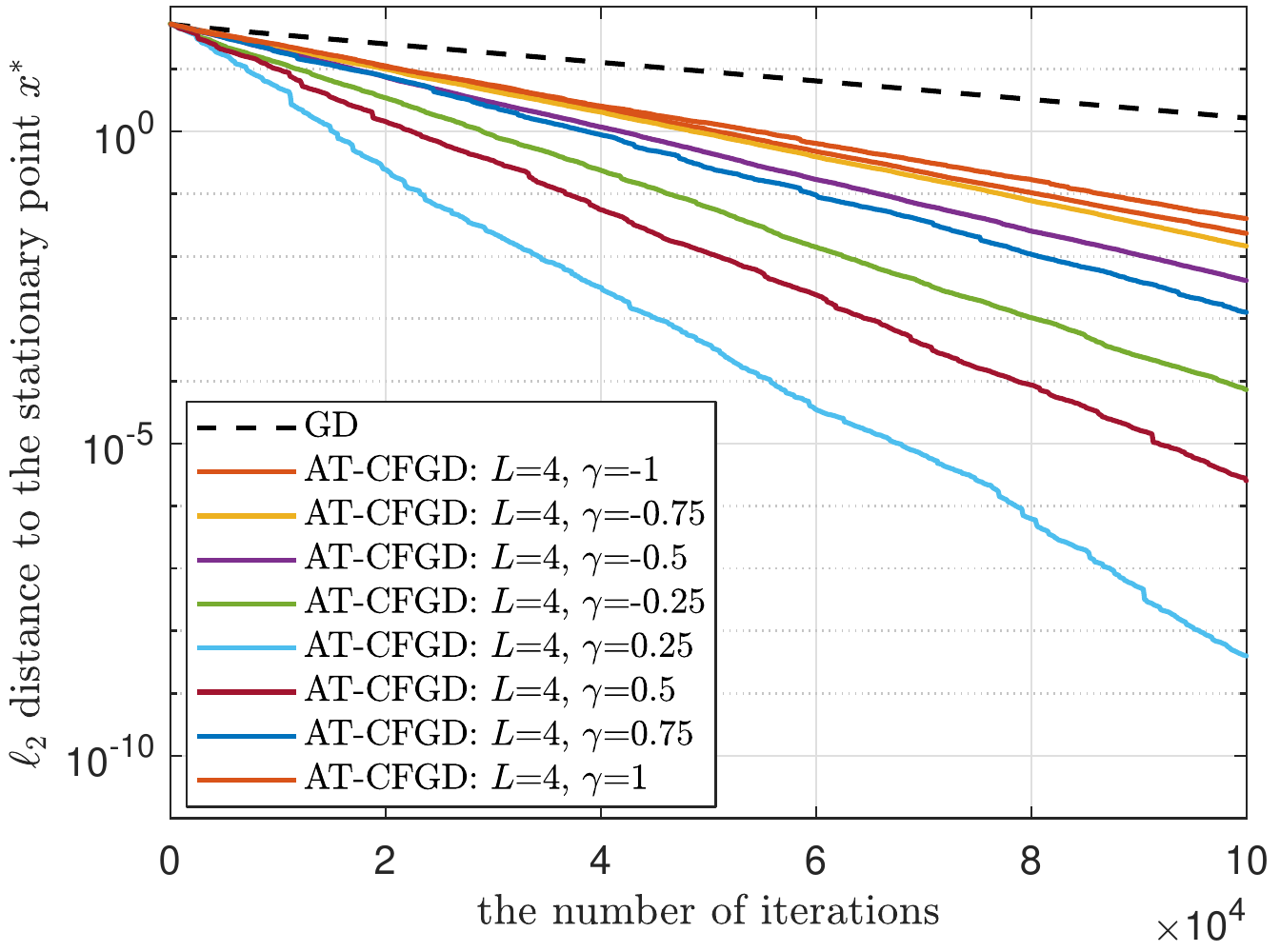}
	}
	\caption{(To be viewed in color)  The $L_2$ distance to the stationary point $\emph{\x}^*$ with respect to the number of iterations at varying $\gamma$ and $L$.
	(Top) $L=1$ and (Bottom) $L=2, 3, 4$ from left to right.
	Gradient descent corresponds to the case of $\gamma = 0$.
	In all cases, the optimal stepsize \eqref{def:LineSearch} from the line search is employed. 
	The condition number of the model matrix $A$ is 90,053.
	}
	\label{fig:convergence-10xm}
\end{figure}

In Figure~\ref{fig:convergence-10xm}, we show the $\ell_2$ distance to the stationary point $\x^*$ 
with respect to the number of iterations
at varying parameters ($L, \gamma$) of AT-CFGD.
Specifically, we report the results of $\gamma \in \{\pm1, \pm0.75, \pm0.5, \pm0.25, 0\}$.
Again, gradient descent corresponds to the special case of $\gamma = 0$.
On the top, the results for $L=1$ are shown. 
We found that AT-CFGD with negative $\gamma$ 
outperforms those with positive $\gamma$ and GD ($\gamma=0$).
In particular,  
AT-CFGD with $\gamma=-0.25$ performs the best among others reaching to the stationary point
within the machine accuracy
in $4\times 10^4$ iterations.
GD exhibits a linear convergence (as expected),
however, since the model matrix $A$ is ill-conditioned (the condition number of $A$ is 90,053), 
the convergence speed is too slow to see 
a significant improvement within $10^5$ iterations.
On the bottom, the results for $L=2,3,4$ are reported.
Unlike the case of $L=1$, 
for all values of $\gamma$, AT-CFGD significantly outperforms GD.
This demonstrates 
the effectiveness of AT-CFGD in mitigating the dependence on the condition number in the rate of convergence. 

\subsection{Quadratic Objective Function: Real data}
We employ the dataset from UCI Machine Learning Repository's ``Gas Sensor Array
Drift at Different Concentrations" \cite{Vergara2012Chemical,Rodriguez2014Calibration}. Specifically, we used the datasets Ethanol
problem a scalar regression task with 2565 examples, each comprising 128 features (one of the largest numeric regression tasks in the repository). The input and
output data sets are normalized to have zero mean and unit variance. After the
normalization, the condition number of the input data matrix $W$ is 70,980, yielding a gigantic condition number of the model matrix $A = WW^\top$ of $(70,980)^2 \approx 5\times 10^{9}$.

\begin{figure}[!htbp]
	\centerline{
	    \includegraphics[width=4.3cm]{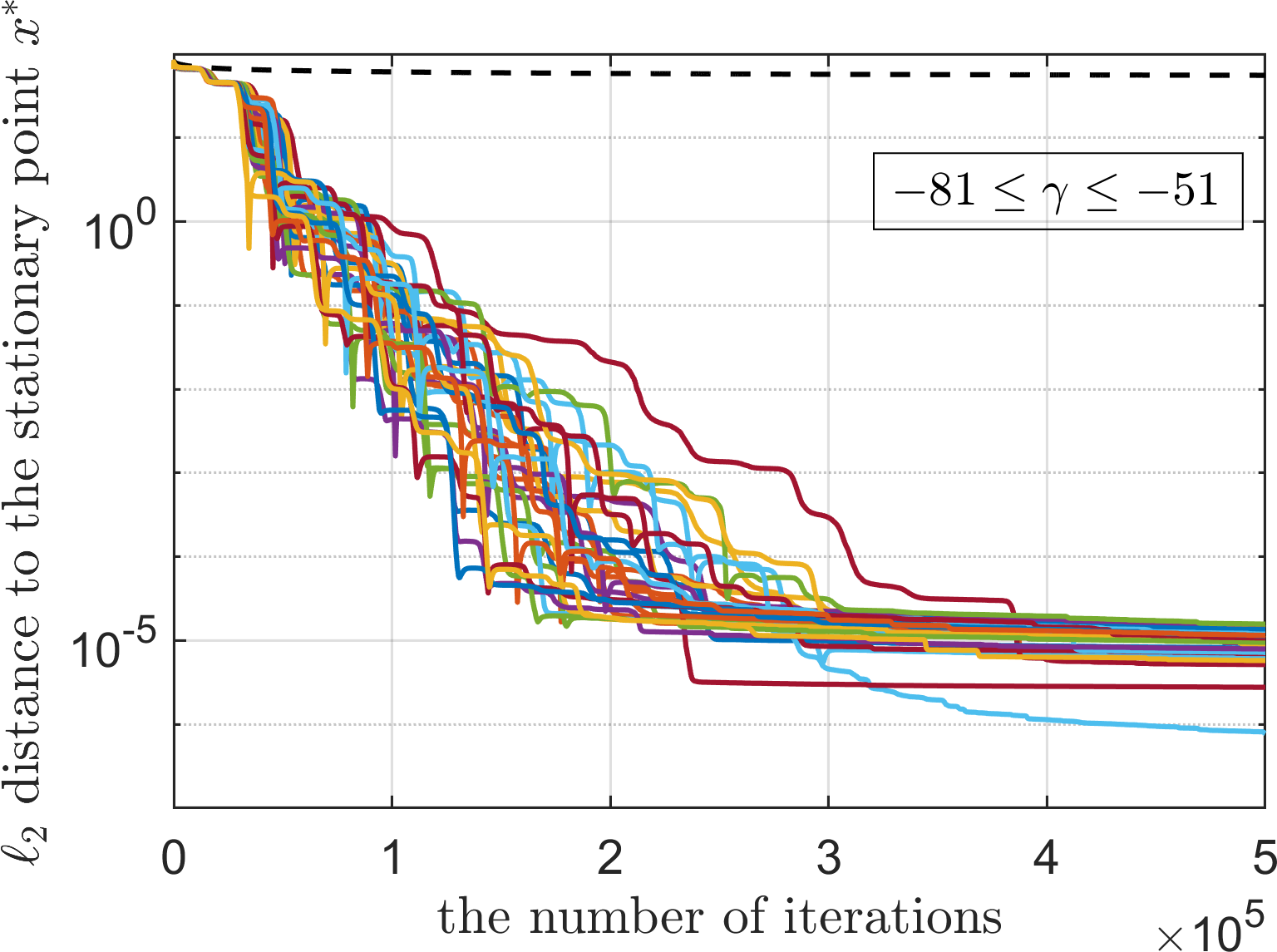}
		\includegraphics[width=4.3cm]{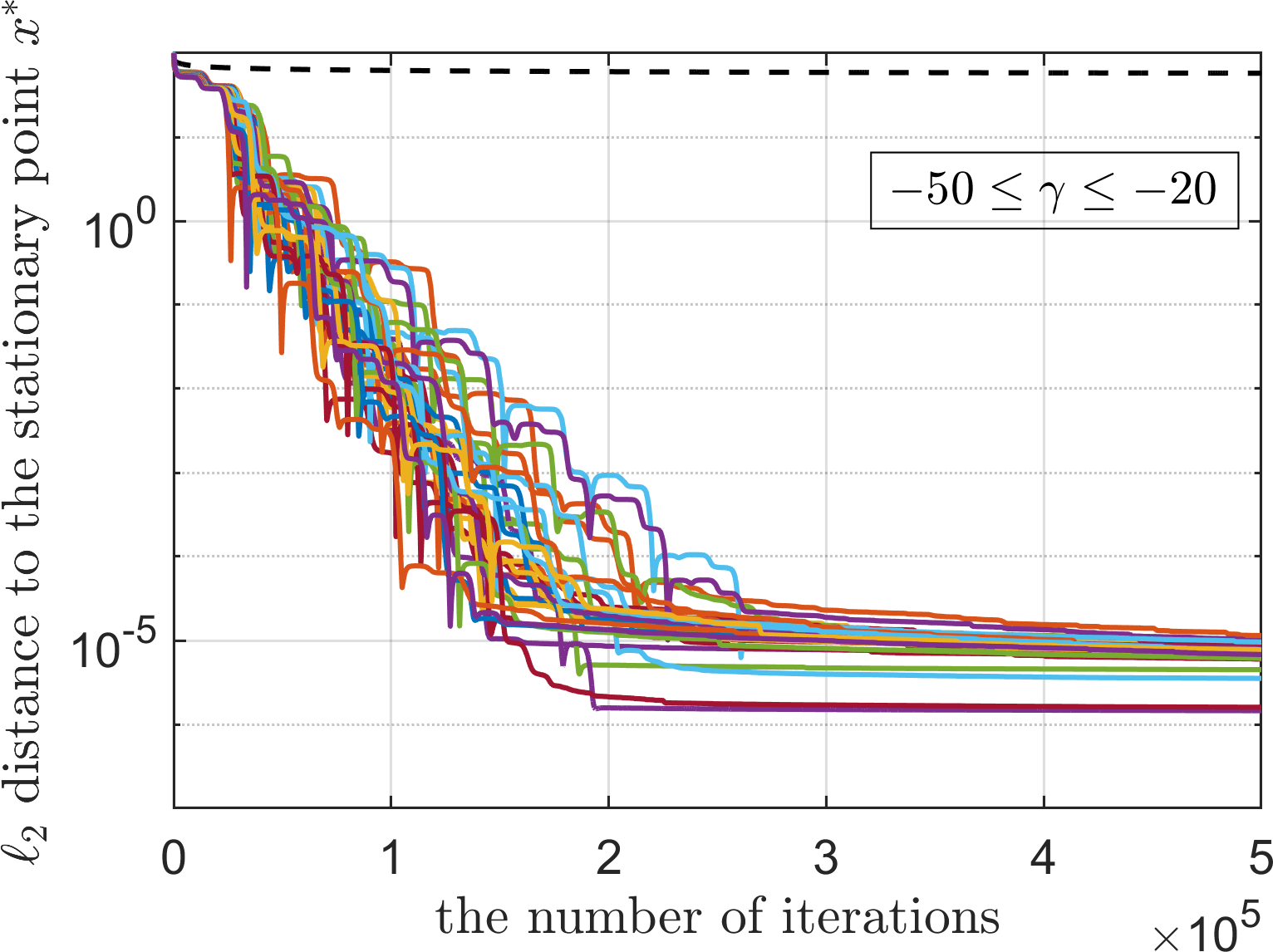}
		\includegraphics[width=4.3cm]{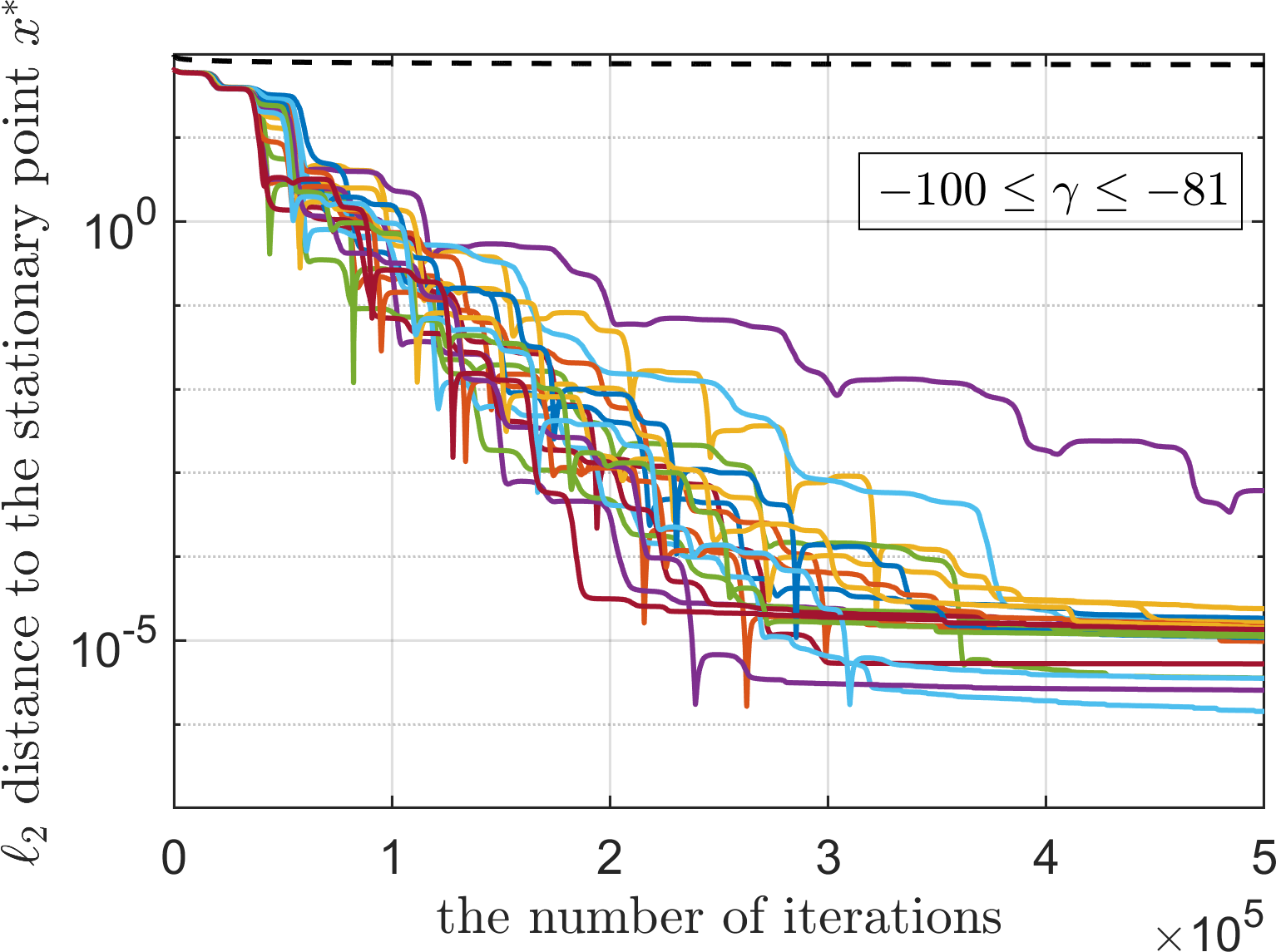}
	}
	\caption{(To be viewed in color) 
	The regression results for the UCI Machine Learning Repository’s dataset of 2565 examples,
	which leads to a gigantic the condition number of $(70,980)^2 \approx 5\times 10^{9}$.
	The $\ell_2$ distance to the stationary point $\emph{\x}^*$ with respect to the number of iterations at varying $\gamma$ and $L=1$.
	Left: $\{\gamma : -80 \le \gamma \le -51\}$.
	Middle: $\{\gamma : -50 \le \gamma \le -20\}$.
	Right: $\{\gamma : -100 \le \gamma \le -81\}$.
	Gradient descent corresponds to the case of $\gamma = 0$ and its trajectories are shown as dashed-lines.
	Both GD and AT-CFGD employ the optimal stepsize \eqref{def:LineSearch}.
	}
	\label{fig:UCI}
\end{figure}

In Figure~\ref{fig:UCI}, we report the $\ell_2$ distance to the stationary point $\x^*$ versus the number of iterations.
We employ AT-CFGD with $L=1$ at varying $\gamma \in [-100, -20]$.
The initial point $\x^{(-1)}$ is set to zero
and $\x^{(0)}$ is randomly chosen from the uniform distribution in the hypercube $[-10,10]^d$.
Gradient descent corresponds to the case of $\gamma = 0$ and its trajectories are shown as black dashed-lines.	
We clearly see that AT-CFGD converges significantly faster than GD.
In almost all cases of $\gamma$,
AT-CFGD converges to the optimum $\x^*$
within the $\ell_2$ error of $10^{-5}$ in merely $3\times 10^{5}$ iterations.
For GD, as expected, we cannot see any significant improvement within $5\times 10^{5}$ iterations. 
Again, this clearly demonstrates that 
AT-CFGD can effectively mitigate the dependence on the condition number in the rate of convergence
and result in a significant acceleration over GD.

\subsection{Nonconvex Objective Functions: Neural Networks} \label{subsec:NNs}

\begin{figure}[!htbp]
	\centerline{
		\includegraphics[width=4.3cm]{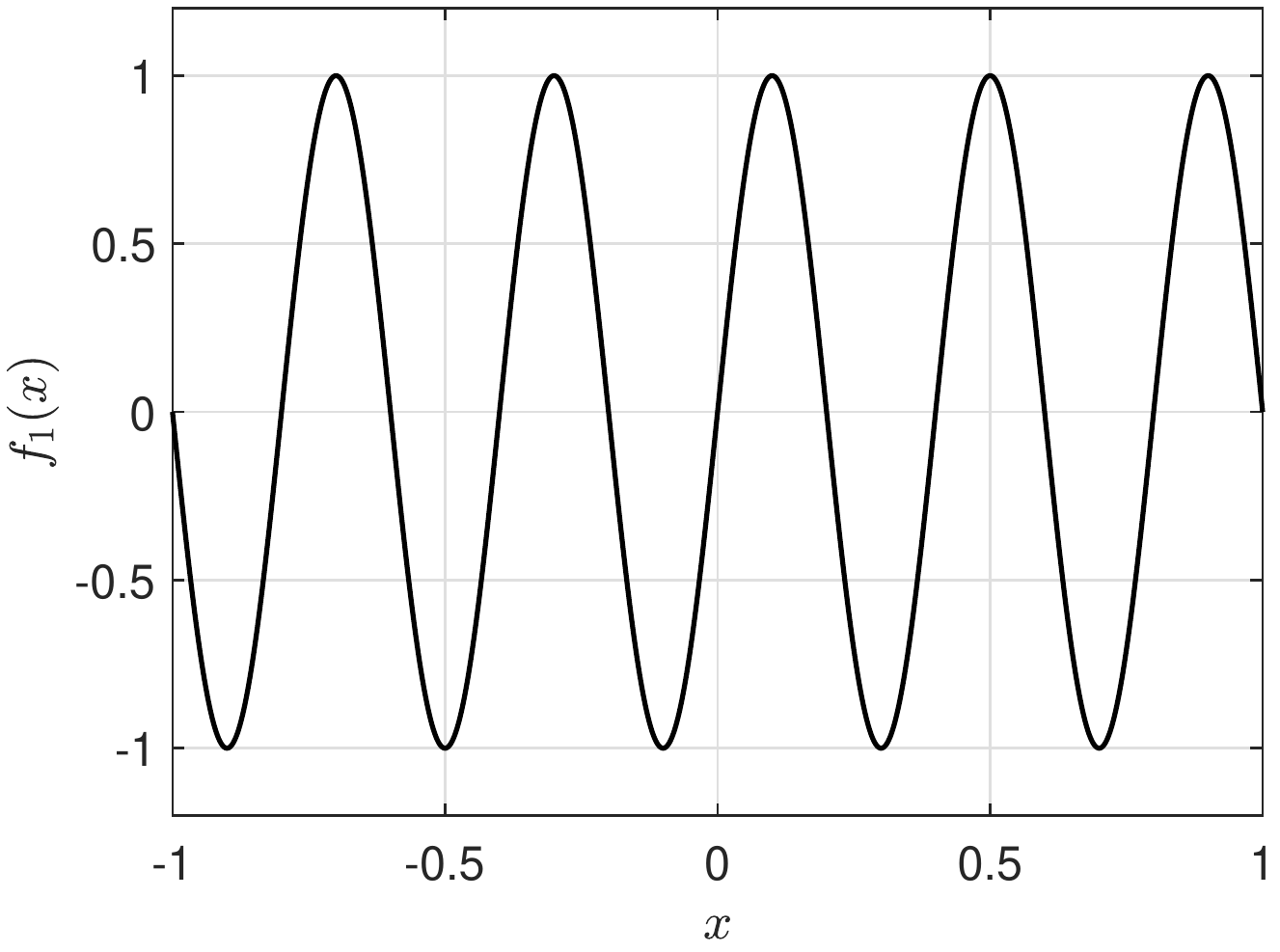}
		\includegraphics[width=4.3cm]{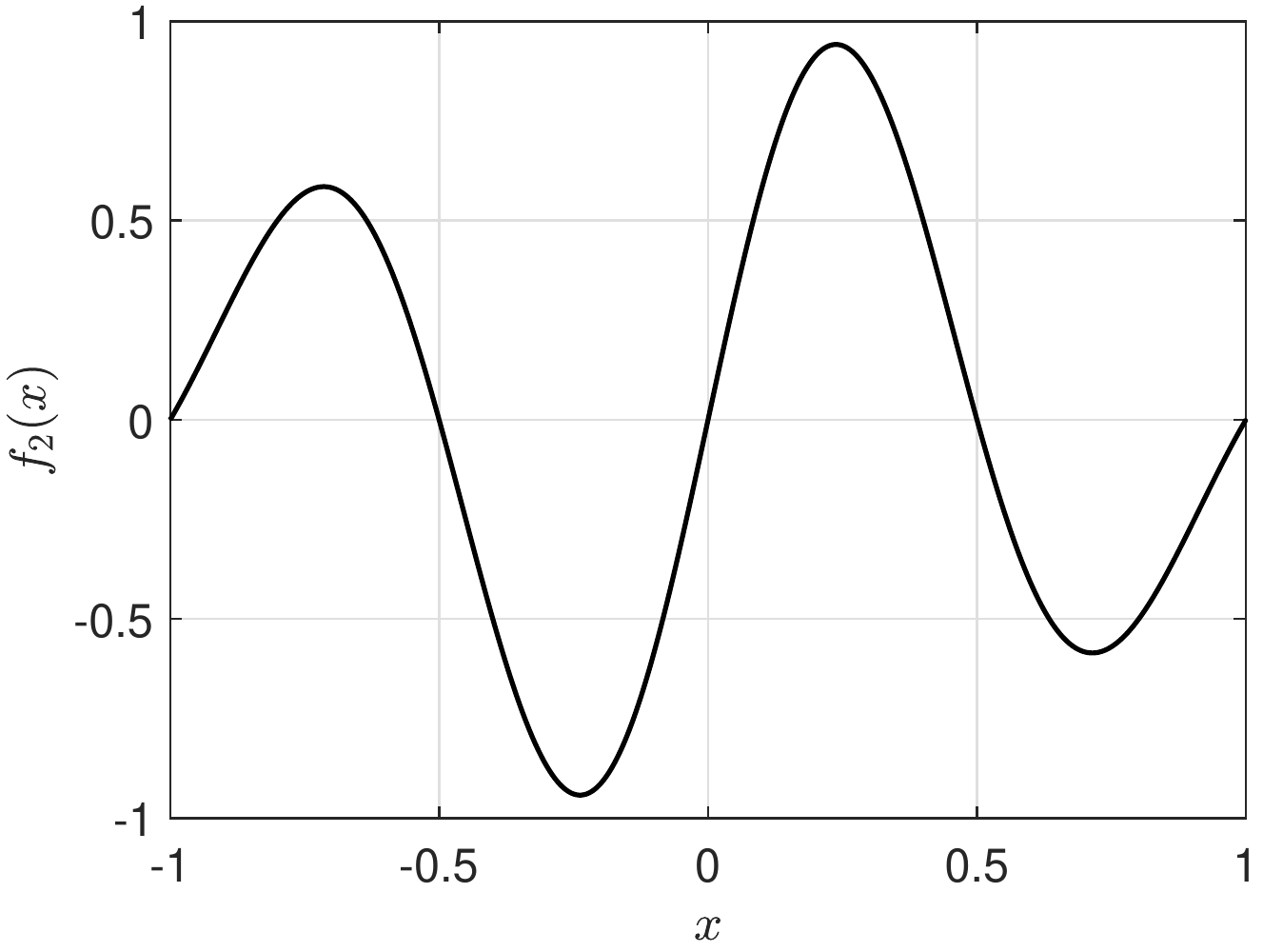}
		\includegraphics[width=4.3cm]{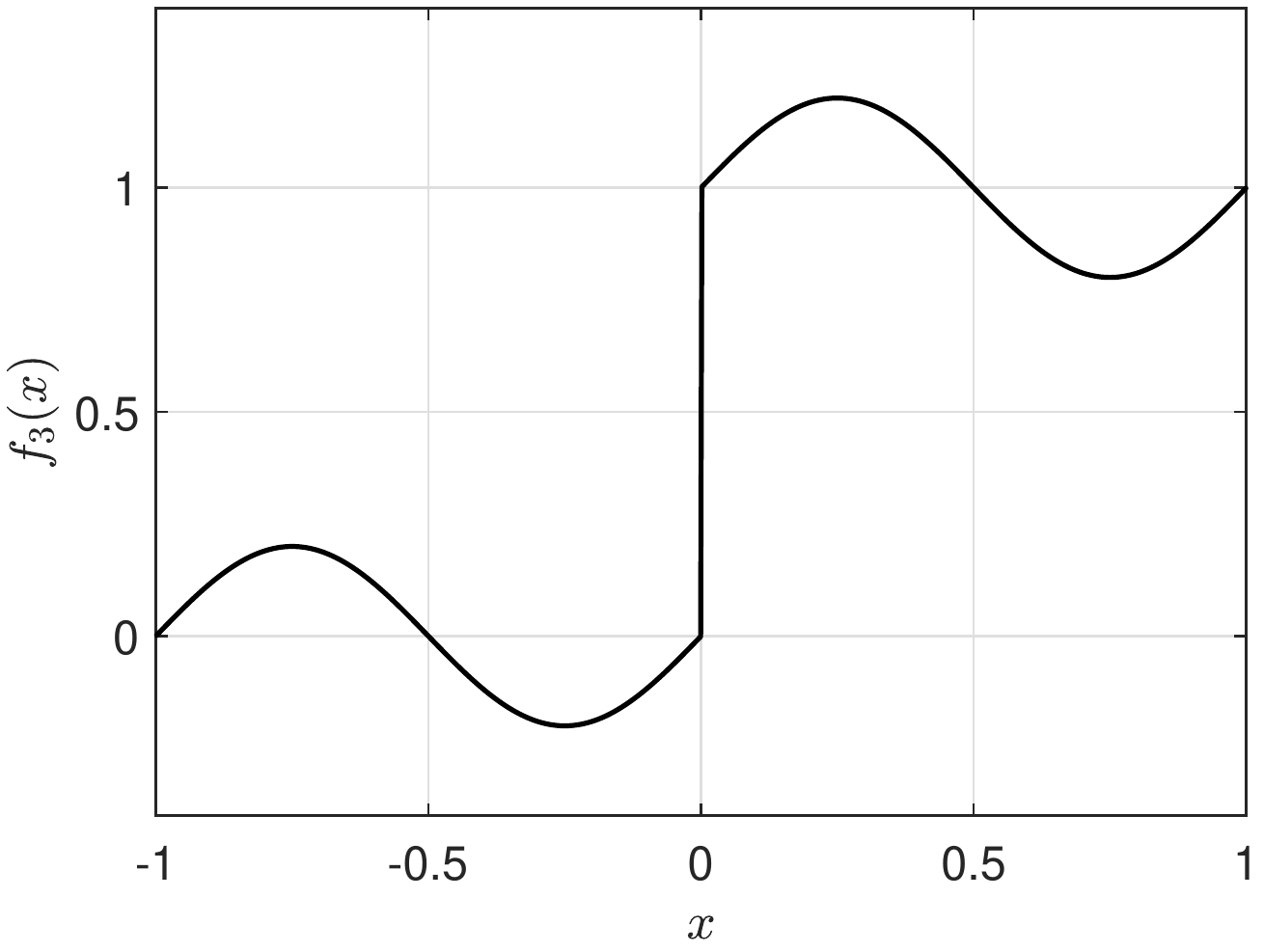}
	}
	\centerline{
		\includegraphics[width=4.3cm]{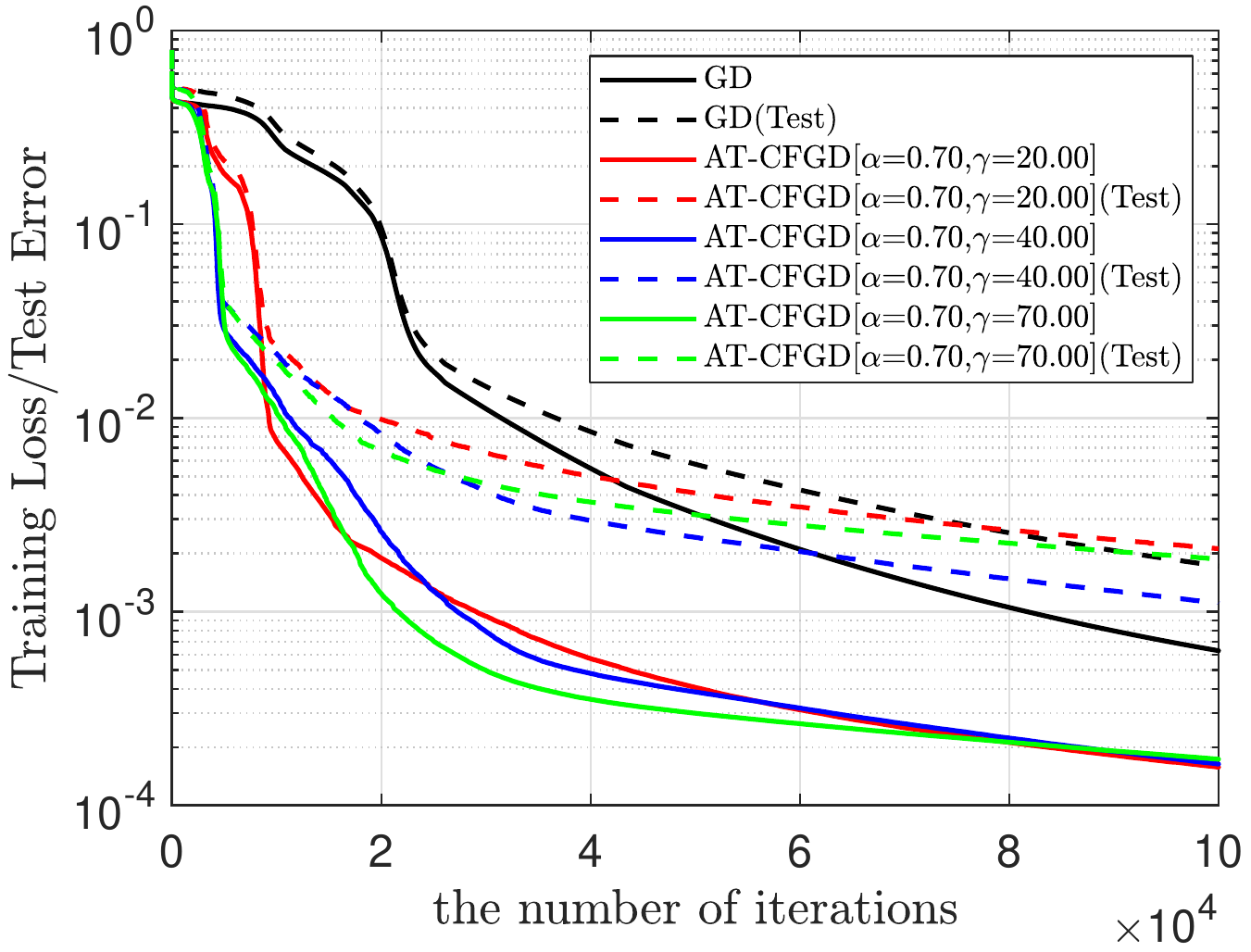}
		\includegraphics[width=4.3cm]{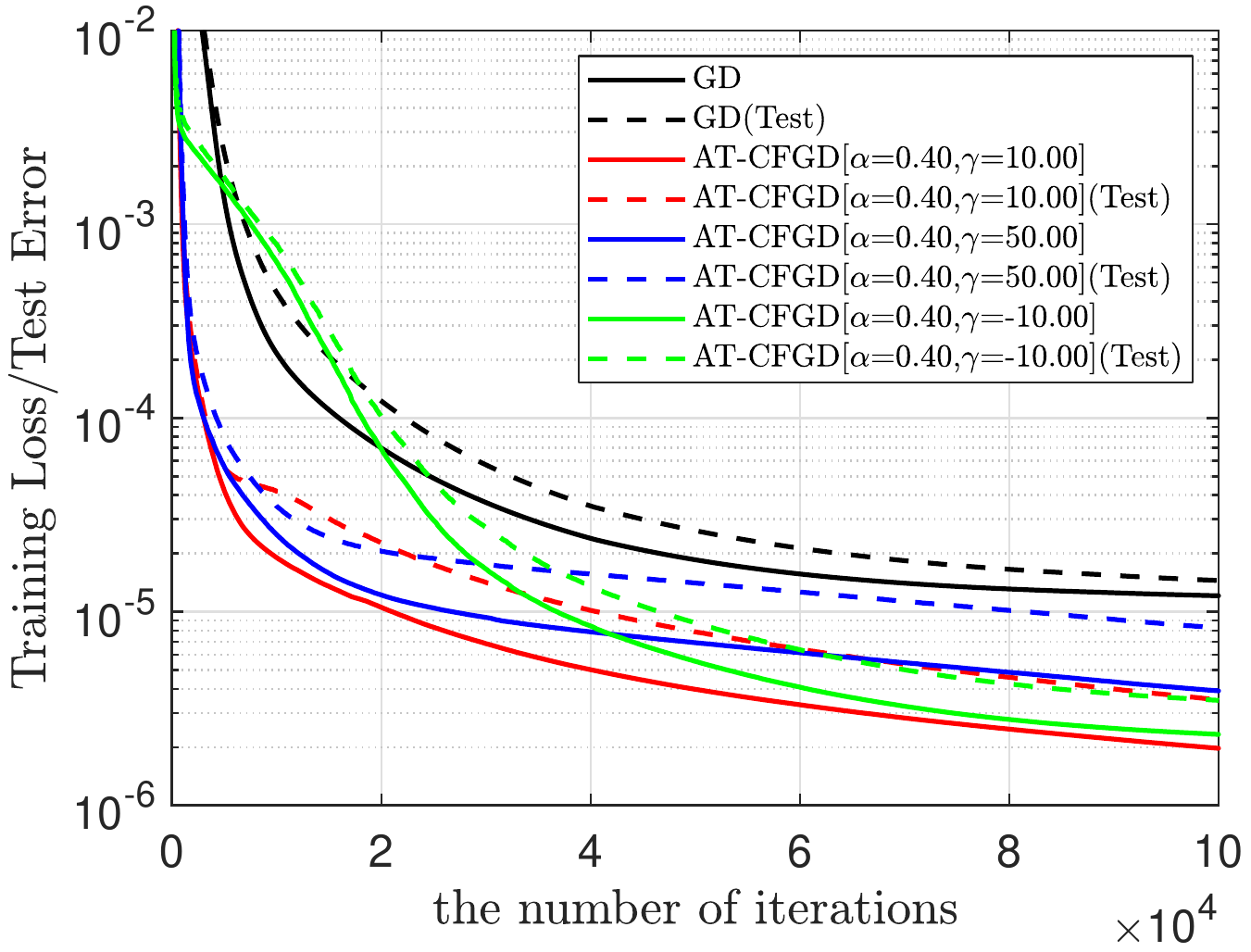}
		\includegraphics[width=4.3cm]{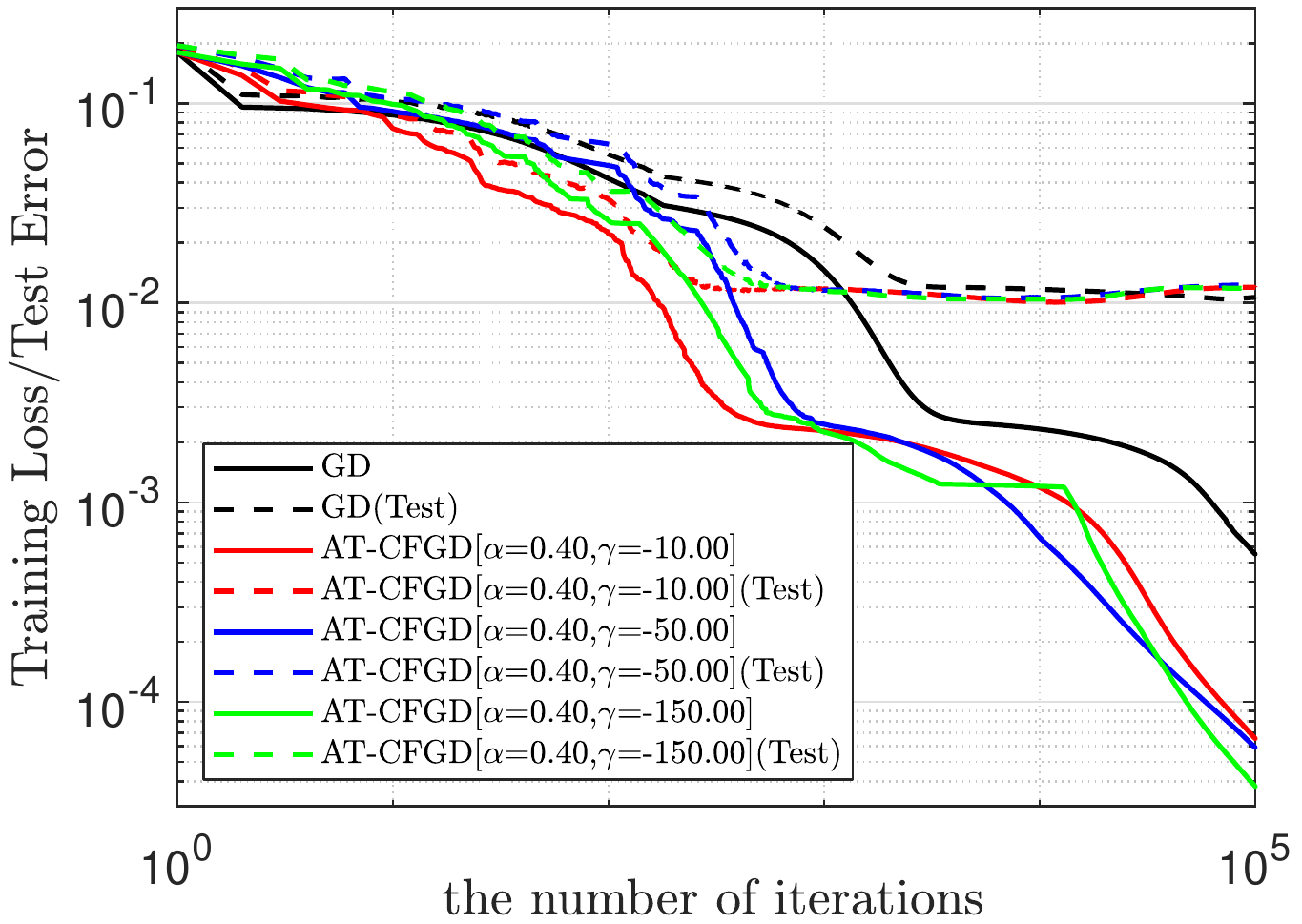}
	}
	\caption{(To be viewed in color) 
	Approximation results for the three test functions \eqref{test-func} -- (left) $h_1(x) = \sin(5\pi x)$, (middle) $h_2(x) = \sin(2\pi x)e^{-x^2}$,
	and (right) $h_3(x) = \mathbb{I}_{x>0}(x) + 0.2\sin(2\pi x)$.
	(Top) The graphs of the three test functions.
	(Bottom) 
	The training loss and the test error trajectories are shown with respect to the number of iterations.
	}
	\label{fig:FA-TEST}
\end{figure}

We consider the training of neural network by CFGD
in function approximation tasks.
The test functions are 
\begin{equation} \label{test-func}
    \begin{split}
        h_1(z) &= \sin(5\pi z), \\
        h_2(z) &= \sin(2\pi z)e^{-z^2}, \\
        h_3(z) &= \mathbb{I}_{z > 0}(z) +0.2\sin(2\pi z),
    \end{split}
\end{equation}
whose graphs are shown in the top row of Figure~\ref{fig:FA-TEST}.
We employ 
the univariate hyperbolic tangent two-layer neural network defined as follow:
\begin{align*}
    N(z;\x) = \sum_{j=1}^n a_{3,j} \tanh(a_{1,j}z + a_{2,j}), \qquad
    \x = \{a_{1,j}, a_{2,j}, a_{3,j}\}_{j=1}^n \in \mathbb{R}^d,
\end{align*}
where $d= 3n$.
Given a set of training data $\{(z_i, h(z_i))\}_{i=1}^m$, 
the objective (loss) function is defined by
\begin{equation}
    f(\x) = \frac{1}{2}\sum_{i=1}^m (N(z_i;\x) - h(z_i))^2.
\end{equation}
Thanks to Theorem~\ref{thm:formula-CFGD},
the Caputo fractional gradient of $f$ can be efficiently computed by using the Gauss-Jacobi quadrature \cite{Ralston_01_NA}.

For fixed $\{a_{1,j}, a_{2,j}\}$, the objective function $f$ with respect to the coefficients $\{a_{3,j}\}$
is quadratic. 
Hence, we employ the optimal learning rate \eqref{def:LSQ-problem} for the coefficients $\{a_{3,j}\}$.
For the weights and the biases $\{a_{1,j}, a_{2,j}\}$,
we select the best learning rate among 32 selections -- $\{t \times 10^{-l} : l=1,\dots,8, t=\frac{1}{4}, \frac{2}{4}, \frac{3}{4}, \frac{4}{4} \}$
in every iteration.
Here the best learning rate is the one that yields the largest decrease in the loss.
For comparison, we also report the results of gradient descent (GD).
Similarly, we employ the optimal learning rate \eqref{def:LSQ-problem} for the coefficients and 
the best one among 32 selections for the weights and biases.
We remark that while the considered learning task is 
univariate function approximation, 
the resulting optimization problem is of dimension  $d=3n$.

In the following tests, we set $m=100$, $n=50$ and use 10 quadrature points.
The 100 training data points are randomly uniformly drawn from $(-1,1)$. 
To measure the performance of the trained neural networks, we also report the test error, which is the mean square error on another 100,000 points uniformly randomly drawn from $(-1,1)$. 
AT-CFGD requires two parameters -- $\alpha$ is a fractional order and $\gamma$ is a parameter that determines $\beta$ such that $\beta = \gamma + \frac{1-\alpha}{2-\alpha}$.

Figure~\ref{fig:FA-TEST}
shows the training loss and the test error trajectories
by GD and AT-CFGD with respect to the number of iterations. 
For AT-CFGD, we set $\alpha = 0.7$ for $h_1$, and $\alpha = 0.4$ for $h_2$ and $h_3$.
For each test, we choose three different values of $\gamma$.
On the left, the approximation results for $h_1$ are reported. 
We observe that AT-CFGD ($\alpha=0.70, \gamma=20, 40, 70$) reaches the loss of $2\times 10^{-4}$ at the end of the training, while 
GD landed at the loss level of $6\times 10^{-4}$. 
At around $2\times 10^4$ iterations, we see that AT-CFGD already reduces the loss function to the level of $10^{-3}$,
while GD has not effectively decrease the loss staying at the level of $10^{-1}$.
This clearly shows the faster convergence of CFGD. 
In the middle and right, we show the results for $h_2$ and $h_3$. 
Again, similar behavior is observed. 
We clearly see that AT-CFGD consistently converges faster than GD in all cases.
Furthermore, we found that AT-CFGD not only converges faster but also produces neural networks that generalize well. 
We see that the test errors are generally smaller than or equal to the one by GD.

\begin{figure}[!htbp]
	\centerline{
		\includegraphics[width=4.3cm]{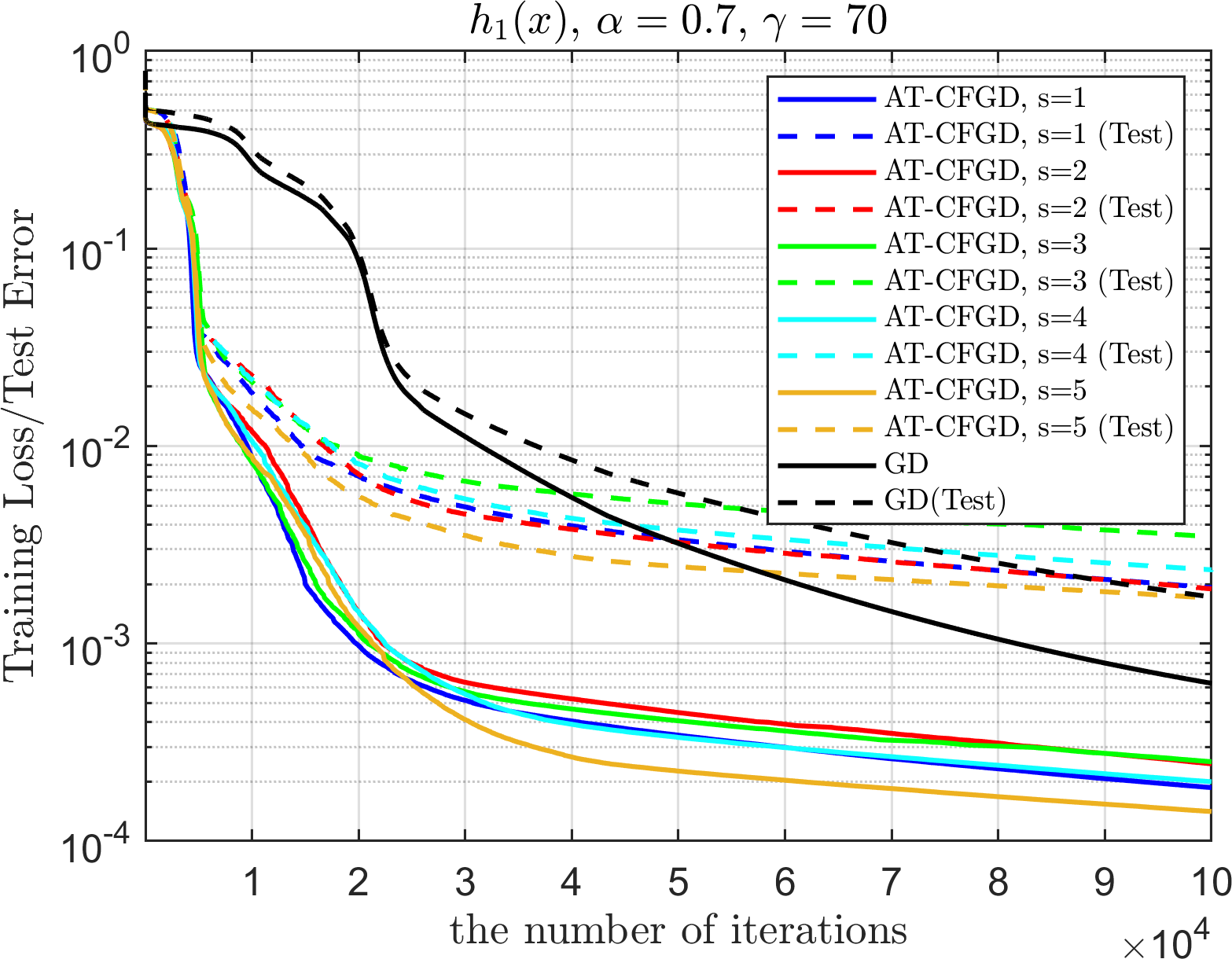}
		\includegraphics[width=4.3cm]{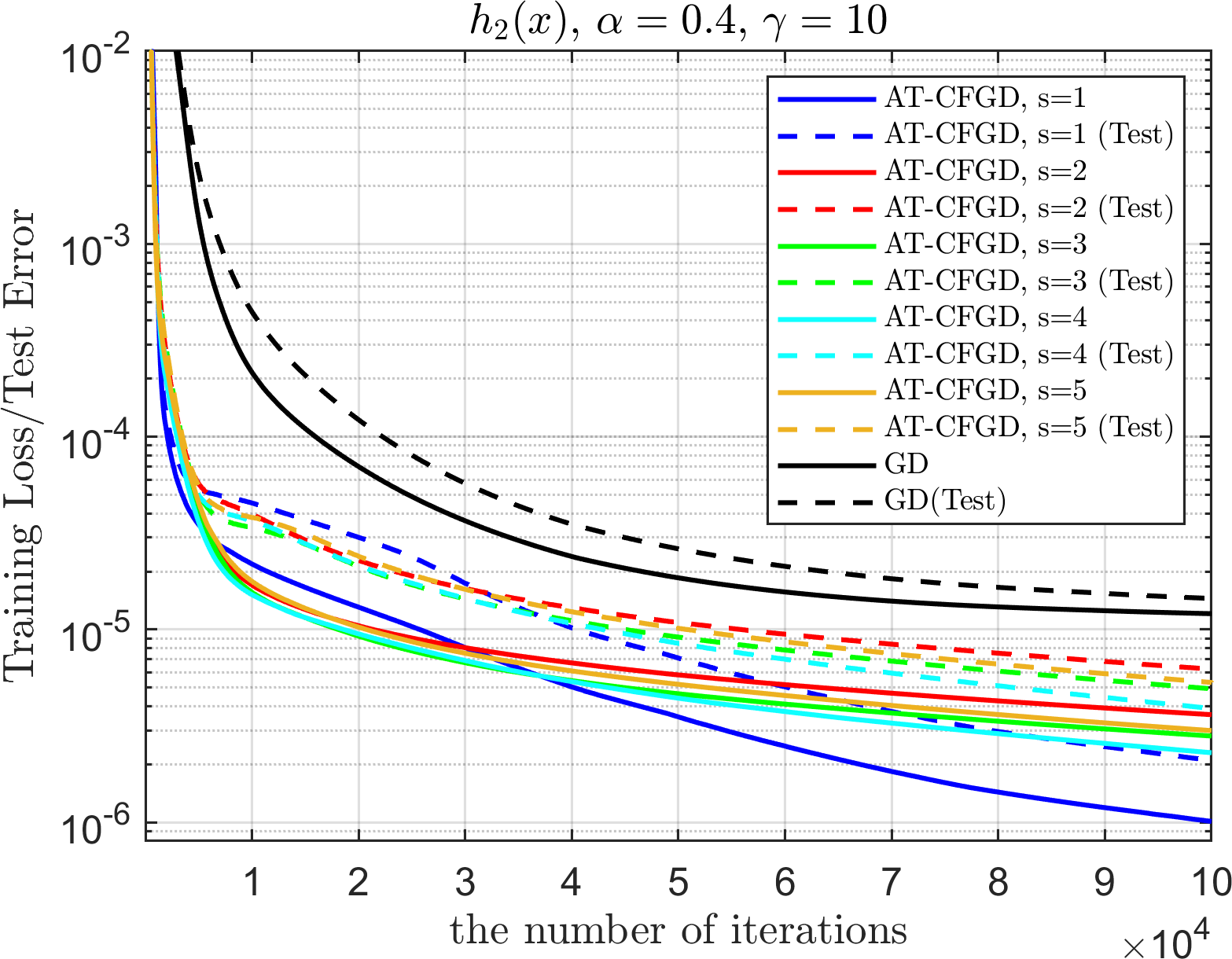}
		\includegraphics[width=4.3cm]{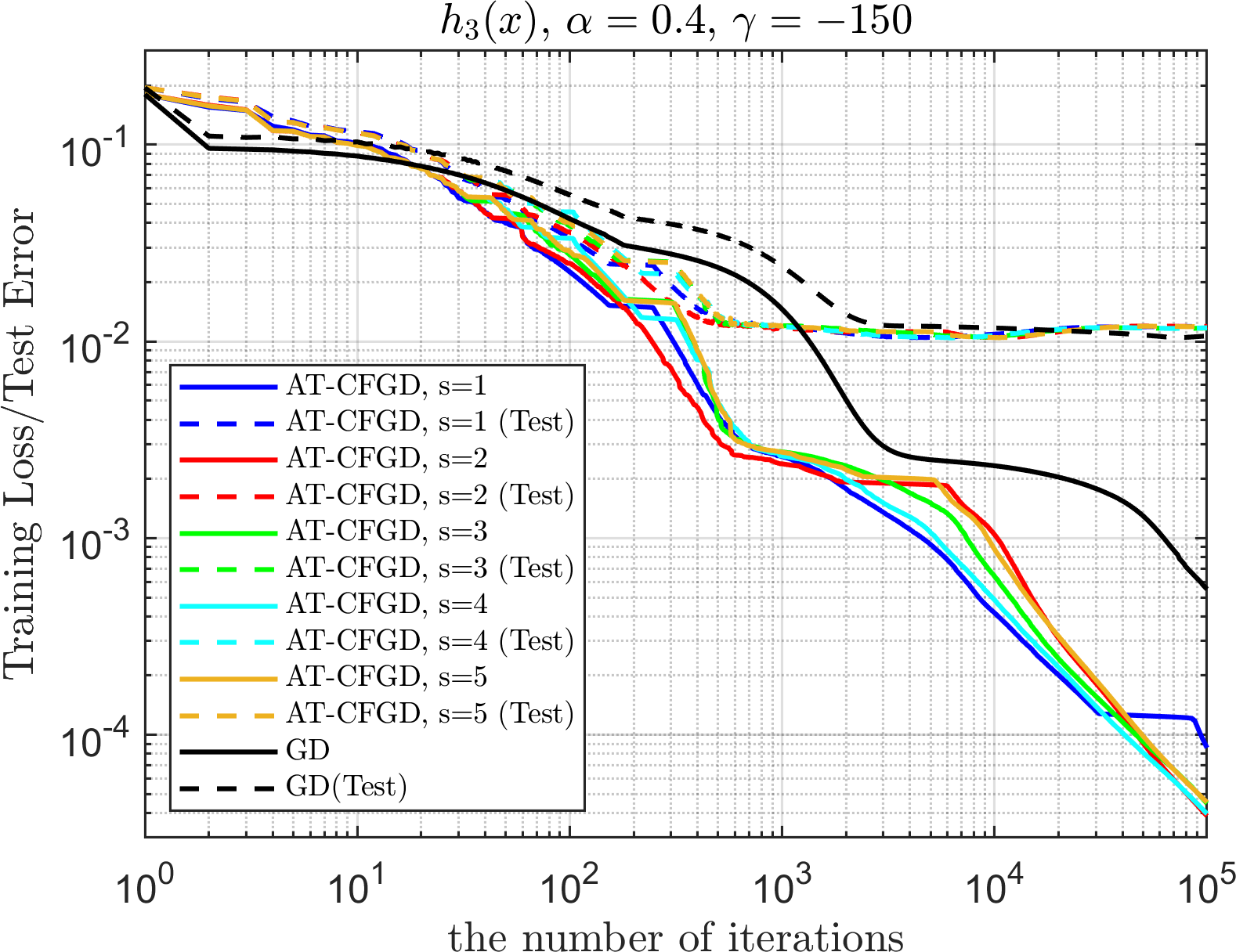}
	}
	\centerline{
		\includegraphics[width=4.3cm]{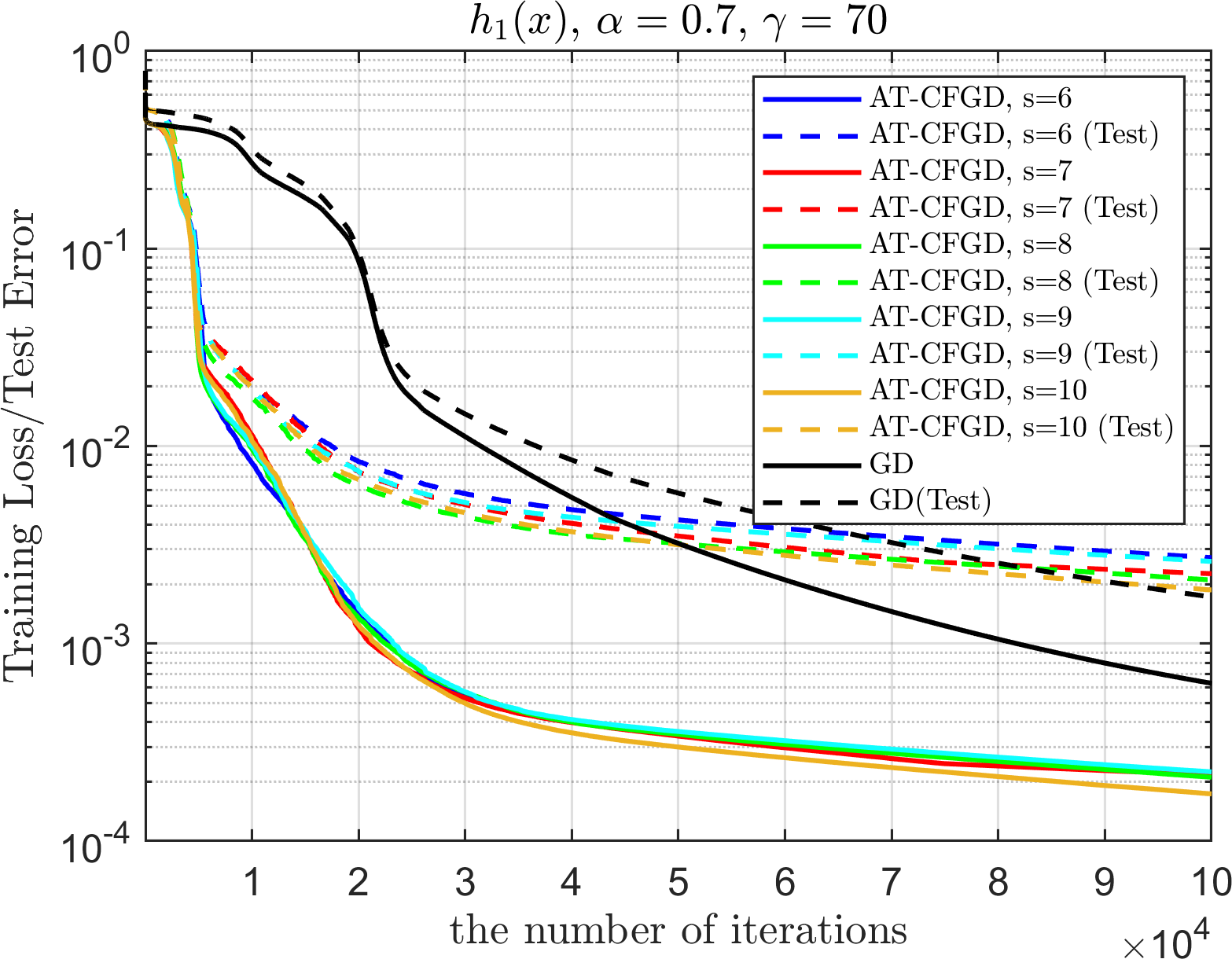}
		\includegraphics[width=4.3cm]{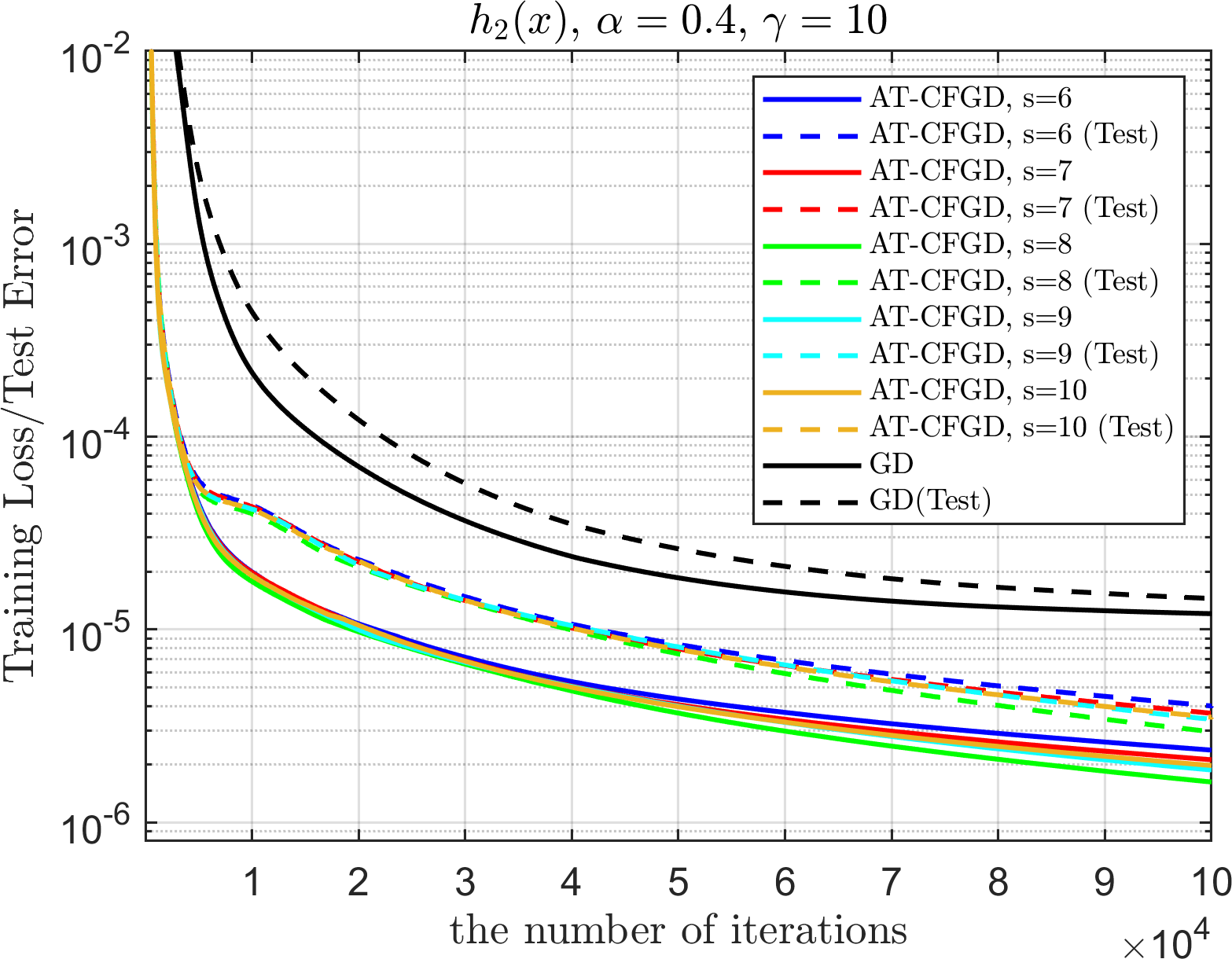}
		\includegraphics[width=4.3cm]{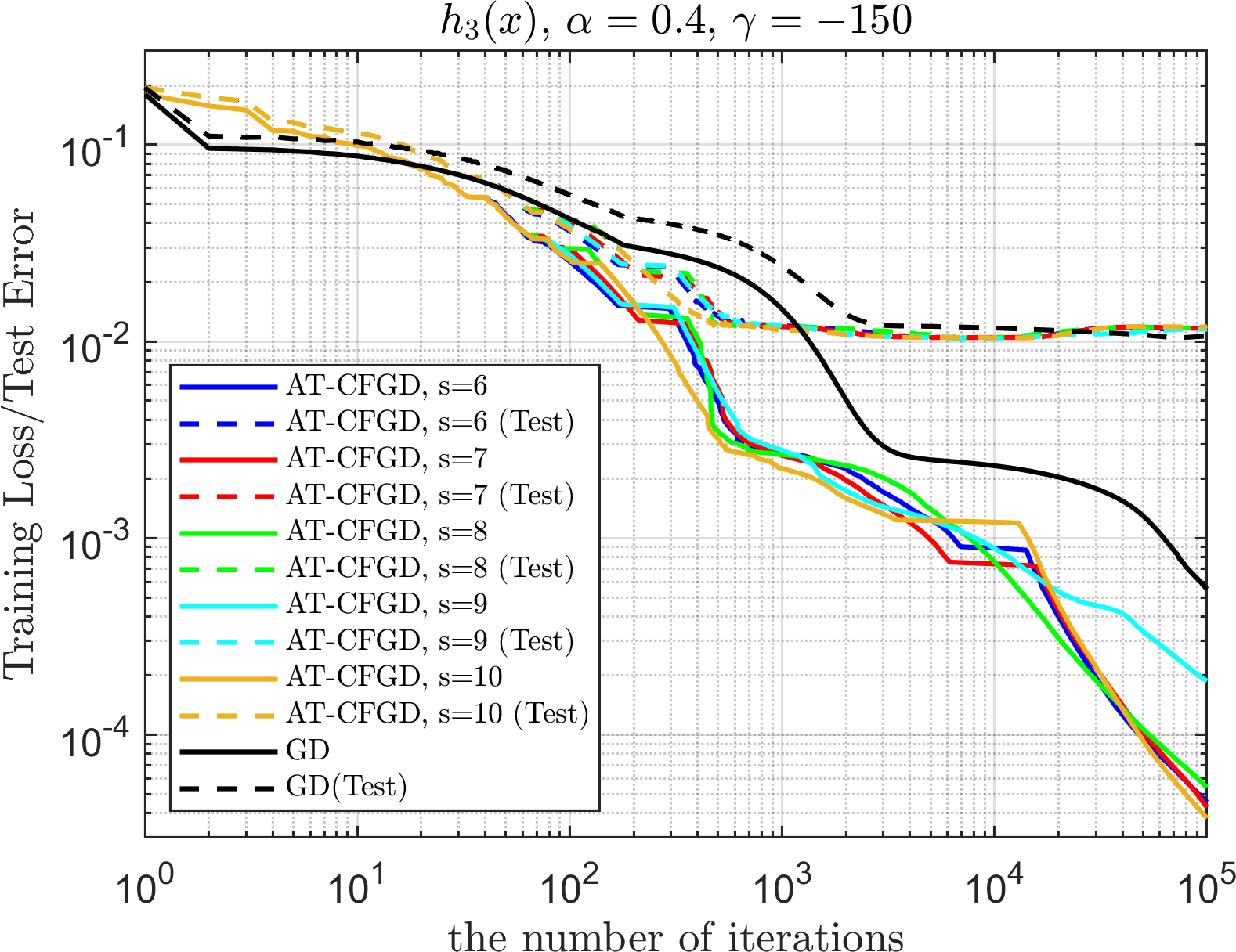}
	}
	\caption{(To be viewed in color) 
	Approximation results for the three test functions \eqref{test-func} -- (left) $h_1(x) = \sin(5\pi x)$, (middle) $h_2(x) = \sin(2\pi x)e^{-x^2}$,
	and (right) $h_3(x) = \mathbb{I}_{x>0}(x) + 0.2\sin(2\pi x)$.
	The training loss and the test error trajectories are shown with respect to the number of iterations
	at varying the number $s$ of the Gauss-Jacobi quadrature points.
	(Top) From $s=1$ to $s=5$.
	(Bottom) From $s=6$ to $s=10$.
	}
	\label{fig:FA-Quad}
\end{figure}

We investigate how the number of quadrature points 
affects the performance of AT-CFGD.
In Figure~\ref{fig:FA-Quad},
we report the results of AT-CFGD implemented by 
using $s$ quadrature points with $s=1,\dots, 10$.
For $h_1$, we set $\alpha = 0.70, \gamma = 70$.
For $h_2$ and $h_3$,
we set $\alpha=0.4$ for both and $\gamma = 10$ and $\gamma=-150$,
respectively. 
We clearly see that for any choices of $s$ including $s=1$,
AT-CFGD outperforms GD.
As discussed in Section~\ref{sec:analysis}, 
the computational cost of AT-CFGD 
is approximately $s$-times higher than those of GD.
Hence, even when the computational cost is taken into account (e.g. $s=1$), 
we see that AT-CFGD converges faster than GD.
We note that when $s$ is small, 
the direction $\vec{\textbf{d}}_k^{\text{QUAD}}$ \eqref{eqn:Quad-CFGD} obtained by a quadrature rule
may no longer be an accurate approximation to 
$\vec{\textbf{d}}_k$ \eqref{eqn:CFGD}.
Yet, we empirically found that regardless of the number of quadrature points, AT-CFGD implemented by $\vec{\textbf{d}}_k^{\text{QUAD}}$ still produces
good directions for the purpose of minimizing the loss function.
We defer further investigation to future work.

Lastly, we found that while AT-CFGD seems not very sensitive to the choice of the parameter $\gamma$,
there is a range of $\gamma$ that makes AT-CFGD more effective.  
Understanding the effect of $\gamma$
also requires 
further investigation, which is deferred to future research.
	
	\appendix
	\section{Basic Calculations}
Here we collected all the necessary calculations involving the Captuo fractional derivative.

\begin{lemma} \label{lem:CFD-x-x2}
    Let $I:\mathbb{R}\to \mathbb{R}$ 
    be the identity map defined by $I(x) = x$.
	For $0 < \alpha < 1$ and for any $c$,
	\begin{align*}
		~{}_{c}^C\!D_{x}^{\alpha} I(x) =\frac{\emph{\text{sign}}(x-c)}{\Gamma(2-\alpha)}|x-c|^{1-\alpha},
		\qquad
		~{}_{c}^{C}\!D_{x}^{\alpha} I^2(x) = 2({}_{c}^{C}D_{x}^{\alpha} x) (\gamma_{\alpha}(x-c) + x),
	\end{align*}
	where $\gamma_{\alpha} = -\frac{1-\alpha}{2-\alpha}$.
	Also, we have
	\begin{align*}
		{}_{c}^{C}D_{x}^{1+\alpha} I(x) &=0,
		\qquad
		{}_{c}^{C}D_{x}^{1+\alpha} I^2(x) 
		= 2~\emph{\text{sign}}(x-c) ~{}_{c}^C\!D_{x}^{\alpha} I(x).
	\end{align*}
\end{lemma}
\begin{proof}
    Direct calculations lead to the results.
\end{proof}

\begin{proposition} \label{prop:caputo-taylor}
	Suppose $f$ is in $C^{\infty}(\mathbb{R})$.
	For $0 < \alpha < 1$ and $c \in \mathbb{R}$,
	we have 
	\begin{equation*}
	\begin{split}
	~{}_{c}^{C}\!D^{\alpha}_{x} f
	&= ({}_{c}^{C}\!D^{\alpha}_{x} I)
	\sum_{k=1}^{\infty} \frac{\Gamma(2-\alpha)}{\Gamma(k+1-\alpha)}f^{(k)}(c)(x-c)^{k-1}, \\
	|x-c|~{}_{c}^{C}\!D^{1+\alpha}_{x} f 
	&= ({}_{c}^{C}\!D^{\alpha}_{x} I)
	\sum_{k=2}^{\infty} \frac{\Gamma(2-\alpha)}{\Gamma(k-\alpha)}f^{(k)}(c)(x-c)^{k-1}.
	\end{split}
	\end{equation*}
\end{proposition}
\begin{proof}
	For $x > c$, by definition and applying integration by parts,
	we obtain 
	\begin{align*}
	~{}_{c}^{C}\!D^{\alpha}_x f
	&=\frac{1}{\Gamma(1-\alpha)}\int_c^x f'(s)(x-s)^{-\alpha} ds \\
	&= 
	\frac{f'(c)}{\Gamma(2-\alpha)}(x-c)^{1-\alpha}
	+ \frac{1}{\Gamma(2-\alpha)}\int_c^x f^{(2)}(s)(x-s)^{1-\alpha} ds.
	\end{align*}
	By repeating integration by parts, we have
	\begin{align*}
	    ~{}_{c}^{C}\!D^{\alpha}_x f
	    &= \sum_{k=1}^{\infty} \frac{f^{(k)}(c)}{\Gamma(k+1-\alpha)}(x-c)^{k-\alpha}
	    \\
	    &= ({}_{c}^{C}\!D^{\alpha}_{x} x) \frac{\Gamma(2-\alpha)}{(x-c)^{1-\alpha}}\sum_{k=1}^{\infty} \frac{f^{(k)}(c)}{\Gamma(k+1-\alpha)}(x-c)^{k-\alpha},
	\end{align*}
	where the second equality uses Lemma~\ref{lem:CFD-x-x2}.
	Suppose $x < c$. Similarly, one can check that	\begin{align*}
	~{}_{x}^{C}\!D^{\alpha}_c f
	&=\frac{-1}{\Gamma(1-\alpha)}\int_x^c f'(s)(s-x)^{-\alpha} ds \\
	&= 
	-\frac{f'(c)}{\Gamma(2-\alpha)}(c-x)^{1-\alpha}
	+ \frac{1}{\Gamma(2-\alpha)}\int_x^c f^{(2)}(s)(s-x)^{1-\alpha} ds
	\\
	&= ({}_{x}^{C}\!D^{\alpha}_c x)f'(c)
	+ \frac{f''(c)}{\Gamma(3-\alpha)}(c-x)^{2-\alpha}
	- \frac{1}{\Gamma(3-\alpha)}\int_x^c f^{(2)}(s)(s-x)^{2-\alpha} ds 
	\\
	&=\cdots 
	=  ({}_{x}^{C}\!D^{\alpha}_c x)\sum_{k=1}^{\infty} \frac{\Gamma(2-\alpha)}{\Gamma(k+1-\alpha)}f^{(k)}(x-c)^k,
	\end{align*}
	which gives one of the desired equations.
	
	The other equation can be checked similarly.
\end{proof}

\section{Proof of Theorem~\ref{thm:smoothing}} \label{app:thm:smoothing}
\begin{proof}
    Let $f$ admit a Taylor expansion at $c$, i.e.,
    \begin{equation*}
    f(z) = \sum_{k=0}^{\infty} \frac{f^{(k)}(c)}{\Gamma(k+1)}(z-c)^k,
    \quad
    f'(z) = \sum_{k=1}^{\infty}         \frac{f^{(k)}(c)}{\Gamma(k)}(z-c)^{k-1}.
    \end{equation*}
    For $\alpha \in (0,1)$, $\beta \in [0,\infty)$ and $c \in \mathbb{R}$, 
    let ${}_{c}F_{\alpha,\beta}$ be a smoothing of $f$
    defined by
    \begin{equation*}
        {}_{c}F_{\alpha,\beta}(z) = f(c) + f'(c)(z-c) + \sum_{k=2}^{\infty} 
        C_{k,\alpha,\beta} \frac{f^{(k)}(c)}{k!}(z-c)^{k},
    \end{equation*}
    where $C_{k,\alpha,\beta} = \left(\frac{\Gamma(2-\alpha)\Gamma(k)}{\Gamma(k+1-\alpha)} + \beta\frac{\Gamma(2-\alpha)\Gamma(k)}{\Gamma(k-\alpha)}
        \right)$.
    The linear approximation of     ${}_{c}F_{\alpha,\beta}(z)$ at $x \ne c$
    is given by $~{}_{c}\ell_{\alpha,\beta}(z) = {}_{c}F_{\alpha,\beta}(x) + {}_{c}F_{\alpha,\beta}'(x) (z-x)$.
    Hence, the steepest descent direction of ${}_{c}F_{\alpha,\beta}(z)$ at $x \ne c$
    is
    $\vec{\textbf{d}} =  {}_{c}F_{\alpha,\beta}'(x)
        = \frac{{}_{c}^{C}\!D^{\alpha}_s f}{{}_{c}^{C}\!D^{\alpha}_s I}
        + \beta |x-c|\frac{{}_{c}^{C}\!D^{1+\alpha}_s f}{{}_{c}^{C}\!D^{\alpha}_s I}$,
    where the second equality holds 
    from Proposition~\ref{prop:caputo-taylor}.
\end{proof}

\section{Proof of Theorem~\ref{thm:formula-CFGD}}
\label{app:thm:formula-CFGD}

The proof is readily followed by the following Lemma.
\begin{lemma} \label{lem:CFD-g}
	Let $f(\emph{\x})$ be
	a real-valued $C^2$ function defined on $\R^d$.
	Let $\alpha \in (0,1)$, $\emph{\textbf{c}} = (c_j)$, $\emph{\x} = (x_j)$
	and $\Delta_j = \frac{|x_j - c_j|}{2}$.
	Then, 
	for $j=1,\dots, d$, we have 
	\begin{align*}
		({}_{c_j}^{C}\!D_{x}^{\alpha} I(x_j))^{-1}
		\big(
		\prescript{C}{\emph{\textbf{c}}}{\nabla_{\emph{\x}}^{\alpha}} f(\emph{\x})\big)_j 
		&= C_{\alpha} 
		\int_{-1}^{1}  f_{j,\emph{\x}}'(\Delta_j(1+u) + c_j)
		(1-u)^{-\alpha} du, \\
		({}_{c_j}^{C}\!D_{x}^{\alpha} I(x_j))^{-1}
		\big(
		\prescript{C}{\emph{\textbf{c}}}{\nabla_{\emph{\x}}^{1+\alpha}} f(\emph{\x})
		\big)_j 
		&= C_{\alpha}
		\int_{-1}^{1}  f_{j,\emph{\x}}''(\Delta_j(1+u) + c_j)
		(1-u)^{-\alpha} du,
	\end{align*}
	where $C_{\alpha} = (1-\alpha)2^{-(1-\alpha)}$
	and $f_{j,\emph{\x}}$'s are defined in \eqref{def:Caputo-gradient}.
\end{lemma}
\begin{proof}
	We will only show the case where $c_j^* < x_j^*$
	as the other case can be done similarly.
	By definition, we have
	$\big(\prescript{C}{\textbf{c}^*}{\nabla_{\x}^{\alpha}} f({\x}^*)\big)_j
		= \frac{1}{\Gamma(1-\alpha)} \int_{c_j^*}^{x_j^*} 
		f_{j,\x^*}'(t) (x_j^* - t)^{-\alpha} dt$.
	Observe that 
	\begin{align*}
		&\frac{1}{\Gamma(1-\alpha)} \int_{c_j^*}^{x_j^*} f_{j,\x^*}'(t) (x_j^* - t)^{-\alpha} dt
		\\
		&= (1-\alpha)\frac{(x_j^*-c_j^*)^{1-\alpha}}{\Gamma(2-\alpha)} 
		\int_{c_j^*}^{x_j^*}  f_{j,\x^*}'(t) (\frac{x_j^* - t}{x_j^*-c_j^*})^{-\alpha} \frac{dt}{x_j^*-c_j^*}
		\\
		&= (1-\alpha){}_{c_j^*}^{C}\!D_{x}^{\alpha} I(x_j^*)
		\int_{0}^1 f_{j,\x^*}'(2\Delta_j s + c_j^*) (1-s)^{-\alpha} ds,
	\end{align*}
	where $\Delta_j = \frac{x_j^* - c_j^*}{2}$, $I(x) = x$ is the identity map
	and 
	the change of variable with $s = \frac{t-c_j^*}{x_j^*-c_j^*}$
	is used in the last equality.
	By further using the change of variable with $u = 2s-1$,
	the above can be written as 
	\begin{align*}
		\frac{1-\alpha}{2^{1-\alpha}}
		~{}_{c_j^*}^{C}\!D_{x}^{\alpha} I(x_j^*)
		\int_{-1}^{1}  f_{j,\x^*}'(\Delta_j(1+u)+c_j^*) 
		(1-u)^{-\alpha} du,
	\end{align*}
	which completes the first part of the proof.
	
	Next, we observe that 
	$\big(
		\prescript{C}{\textbf{c}^*}{\nabla_{\x}^{1+\alpha}} f({\x}^*)
		\big)_j
		= \frac{1}{\Gamma(1-\alpha)} \int_{c_j^*}^{x_j^*} f_{j,\x^*}''(t) (x_j^* - t)^{-\alpha} dt$.
	It then can be checked that 
	\begin{align*}
		&\frac{1}{\Gamma(1-\alpha)} \int_{c_j^*}^{x_j^*} f_{j,\x^*}''(t) (x_j^* - t)^{-\alpha} dt
		\\
		&= (1-\alpha)\frac{(w_j^*-c_j^*)^{1-\alpha}}{\Gamma(2-\alpha)} 
		\int_{c_j^*}^{x_j^*}  f_{j,\x^*}''(t) (\frac{x_j^* - t}{x_j^* - c_j^*})^{-\alpha} \frac{dt}{x_j^*-c_j^*}
		\\
		&=(1-\alpha) ~{}_{c_j^*}^{C}\!D_{x}^{\alpha} I(x_j^*)
		\int_0^1 f_{j,\x^*}''(2\Delta_js + c_j^*)(1-s)^{-\alpha}ds.
	\end{align*}
	By using the change of variable with $u = 2s-1$,
	the above can be written as 
	\begin{align*}
		\frac{1-\alpha}{2^{1-\alpha}}
		~{}_{c_j^*}^{C}\!D_{x}^{\alpha} I(x_j^*)
		\int_{-1}^{1}  f_{j,\x^*}''(\Delta_j(1+u) + c_j^*)(1-u)^{-\alpha}du.
	\end{align*}
\end{proof}

\section{Proof of Theorem~\ref{thm:fLSQ}} \label{app:thm:fLSQ}
\begin{proof}
    Let $\bc = (c_j)$ and $\x = (x_j)$.
    Let 
	\begin{equation*}
	    M = \text{diag}(
	    \begin{bmatrix}
	        {}_{c_1}^C\! D_{x_1}^{\alpha} x_1
	        &
	        \cdots
	        & 
	        {}_{c_d}^C\! D_{x_d}^{\alpha} x_d
	    \end{bmatrix}
	    ) \in \mathbb{R}^{d\times d}.
	\end{equation*}
	For $0 < \alpha < 1$, it can be checked that 
	\begin{align*}
	\big(\prescript{C}{\bc}{\nabla_{\x}^{\alpha} f(\x)}\big)_j
	&=~{}_{c_j}^{C}\!D_{x_j}^{\alpha}x_j
	\sum_{i=1}^m \bigg(w_{ij}^2 \gamma_{\alpha}(x_j - c_j) 
	+  w_{ij}(\w_i^\top \x -y_i) 
	\bigg),
	\end{align*}
	where the equality uses 
	Lemma~\ref{lem:CFD-x-x2}
	and $\gamma_\alpha = -\frac{1-\alpha}{2-\alpha}$.
	Then, the Caputo fractional gradient of $f(\x)$ is given by 
	${}_{\bc}^{C}\nabla_{\x}^{\alpha} f(\x) 
    = M
	\left[W(W^\top \x - \y) 
	+ \gamma_{\alpha} \text{diag}(\tilde{R})(\x-\bc)\right]$.
	Similarly, 
	$\big(\prescript{C}{\bc}{\nabla_{\x}^{1+\alpha} f(\x)}\big)_j
	= \text{sign}(x_j-c_j)~{}_{c_j}^{C}\!D_{x_j}^{\alpha}x_j
	\sum_{i=1}^m w_{ij}^2$,
	which gives 
	${}_{\bc}^{C}\nabla_{\x}^{1+\alpha} f(\x)
	    = M\text{diag}(\text{sign}(x_j-c_j))\tilde{R}$,
	where $\tilde{R} = (\sum_{k=1}^m w_{kj}^2)_j \in \mathbb{R}^d$.
	Then,
	\begin{align*}
	    &{}_{\bc}^{C}\nabla_{\x}^{\alpha} f(\x)
	    +
	    \beta ~\text{diag}(|x_j-c_j|)~{}_{\bc}^{C}\nabla_{\x}^{1+\alpha} f(\x) \\
	    &= 
	    M\left[W(W^\top \x - \y) 
	+ (\beta + \gamma_{\alpha}) \text{diag}(\tilde{R})(\x-\bc)\right].
	\end{align*}
	Let $\gamma_{\alpha,\beta} = \beta + \gamma_{\alpha}$.
	We then have 
	$\vec{\textbf{d}} = 
	WW^\top\left[ 
	(I + \gamma_{\alpha,\beta} K)\x - (\x^* + \gamma_{\alpha,\beta} K \bc)
	\right]$,
	where 
	$K = (WW^\top)^{-1}\text{diag}(\tilde{R})$.
	Let 
	$\x^*_{\alpha,\beta} = (I+ \gamma_{\alpha,\beta} K)^{-1}(\x^* + \gamma_{\alpha,\beta} K \bc)$.
	It then can be checked that 
	$\vec{\textbf{d}} = 
	(WW^\top + \gamma_{\alpha,\beta} \text{diag}(\tilde{R}))
	(\x - \x^*_{\alpha,\beta})$.
	Let $A_{\alpha,\beta} = WW^\top + \gamma_{\alpha,\beta} \text{diag}(\tilde{R})$.
	Note that 
	$[A_{\alpha,\beta} ]_{ij} = (\beta+ \frac{1}{2-\alpha})\sum_{k=1}^m w_{ki}^2$
	if $i = j$
	and $\sum_{k=1}^m w_{ki}w_{kj}$ otherwise.
	It follows from the Caputo fractional gradient descent that
	\begin{align*}
	\x^{(k+1)} - \x^*_{\alpha,\beta} = 
	(I - \eta_k A_{\alpha,\beta})(\x^{(k)} - \x^*_{\alpha,\beta}).
	\end{align*}
	Since $\alpha$ and $\beta$ are chosen to make $A_{\alpha,\beta}$ positive definite,
	let $\sigma_{\max}$ be the largest singular values of $A_{\alpha,\beta}$.
	Let $\kappa$ be the condition number of $A_{\alpha,\beta}$.
	Suppose $\eta_k = \frac{\eta}{\sigma_{\max}}$ 
	for some $\eta \in (0,2)$.
	Then, we have
	$\|\x^{(k)} - \x_{\alpha,\beta}^*\|^2 \le 
	\|\x^{(0)} - \x_{\alpha,\beta}^*\|^2 
	|1- \frac{\eta}{\kappa}|^k$.
	By observing that 
	\begin{align*}
	    \x^*_{\alpha,\beta} &= 
	    \bc + (I+ \gamma_{\alpha,\beta} K)^{-1}(\x^* - \bc)
	    \\
	    &=
	    \bc + (WW^\top+ \gamma_{\alpha,\beta} \text{diag}(\tilde{R}))^{-1}WW^\top(\x^* - \bc) \\
	    &=\bc + (WW^\top+ \gamma_{\alpha,\beta} \text{diag}(\tilde{R}))^{-1}W(\y - W^\top \bc) = \x^*_{\text{Tik}},
	\end{align*}
	the proof is completed.
\end{proof}

\section{Proof of Theorem~\ref{thm:adapt-terminal}} \label{app:thm:adapt-terminal}
\begin{proof}
	For a positive integer $L$, we observe that 
	\begin{align*}
		\vec{\textbf{d}}_k
		= A\x^{(k)} + b + {\gamma}_{\alpha,\beta}\text{diag}(A) (\x^{(k)}-\x^{(k-L)}),
	\end{align*}
	where ${\gamma}_{\alpha,\beta} = \beta - \frac{1-\alpha}{2-\alpha}$
	and $\text{diag}(A)$ is the diagonal matrix whose diagonal entries are from $A$.
	Since $\x^* = -A^{-1}b$, we have
	\begin{align*}
		\x^{(k+1)} - \x^*
		= \left[I - \eta (A + {\gamma}_{\alpha,\beta} \text{diag}(A)) \right](\x^{(k)}-\x^*) 
		+ \eta {\gamma}_{\alpha,\beta}\text{diag}(A) (\x^{(k-L)} - \x^*).
	\end{align*}
	Let $\mathcal{A} = I-\eta (A + \mathcal{B})$,
	$\mathcal{B} =  {\gamma}_{\alpha,\beta}\text{diag}(A)$
	and $\Delta^k = \x^{(k)} -\x^*$.
	Then, the above can be written as $\mathcal{E}_{k+1} = M\mathcal{E}_{k}$ where 
	\begin{align*}
		\mathcal{E}_k=\begin{bmatrix}
			\Delta^{k} &  \cdots & \Delta^{k-L}
		\end{bmatrix}^\top,
	    \quad 
		M = \begin{bmatrix}
			\vec{\mathcal{A}} &  \mathcal{B} \\ 
			I_{Ld} & 0 
		\end{bmatrix},
		\quad
		\vec{\mathcal{A}} = \begin{bmatrix}
		\mathcal{A} & 0_{d \times (L-1)d} 
		\end{bmatrix}.
	\end{align*}
	Here $0_{d \times (L-1)d}$ is the zero matrix of size $d \times (L-1)d$
	and $I_{Ld}$ is the identity matrix of size 
	$Ld$.
	We then obtain  
	\begin{equation} \label{eqn:AT-error}
		\Delta^k = \begin{bmatrix}
	        I & 0_{d\times Ld}
	    \end{bmatrix}
	    \mathcal{E}_k 
	    = \begin{bmatrix}
	        I & 0_{d\times Ld}
	    \end{bmatrix} M^k \mathcal{E}_0.
	\end{equation}
	Let $M^k = \begin{bmatrix} \mathcal{C}^0_k & \mathcal{C}^1_k & \cdots & \mathcal{C}^{L}_k \end{bmatrix}$
	where $\mathcal{C}^{j}_k \in \mathbb{R}^{(L+1)d \times d}$. 
	It then can be checked that 
	\begin{align*}
	    \begin{bmatrix} \mathcal{C}^0_{k+1} & \mathcal{C}^1_{k+1} & \cdots & \mathcal{C}^{L-1}_{k+1} & \mathcal{C}^{L}_{k+1} \end{bmatrix} = 
	    \begin{bmatrix} \mathcal{C}^0_k\mathcal{A} + \mathcal{C}^1_k & \mathcal{C}^2_k & \cdots & \mathcal{C}^{L}_k & \mathcal{C}^0_k\mathcal{B}
	    \end{bmatrix},
	\end{align*}
	which gives the following recurrent relations:
	Let the first $d$ rows of $M^k$ be 
	\begin{align*}
	    \begin{bmatrix}
	        I & 0_{d\times Ld}
	    \end{bmatrix} M^k
	    = \begin{bmatrix}
	\mathcal{A}_k^0 & \mathcal{A}_k^1 & \cdots & \mathcal{A}_k^{L}
	\end{bmatrix}.
	\end{align*} 
	Then, starting with $\mathcal{A}_{1}^0 = \mathcal{A}$, 
	$\mathcal{A}_{1}^L = \mathcal{B}$
	and $\mathcal{A}_{1}^j = 0$ for $1 \le j < L$,
	we have, 
	\begin{align*}
	    \mathcal{A}_k^0 = \mathcal{A}_{k-1}^0\mathcal{A}_1^0 + \mathcal{A}_{k-1}^1, \qquad
	    \mathcal{A}_k^j = \mathcal{A}_{k-1}^{j+1}, \quad \forall 1 \le j < L, \qquad
	    \mathcal{A}_k^{L} = \mathcal{A}_{k-1}^0\mathcal{A}_1^L,
	\end{align*}
	for $k=2,\dots$.
	It then follows from \eqref{eqn:AT-error} that 
	$ \|\Delta^k\| \le \sum_{j=0}^{L} \|\mathcal{A}_k^{j}\|\|\Delta^{-j}\|$.
	
	Since 
    $\mathcal{A}_{k+s,L-s+1} = \mathcal{A}_{k+1,L} = \mathcal{A}_{k,0}\mathcal{A}_{1,L}$ 
    for $s=1,\dots, L$,
	if $\|\mathcal{A}_{k,0}\|$ converges to 0 as $k\to \infty$,
    AT-CFGD converges to the optimal solution to \eqref{def:Quad-problem}
    and the proof is completed.
\end{proof}

\section{Proof of Theorem~\ref{thm:adapt-order}} \label{app:thm:adapt-order}

\begin{proof}
	Let $\x^*_{\gamma}$ be the solution to the Tikhonov regularization \eqref{def:tik-opt}, where
	$\gamma = \beta - \frac{1-\alpha}{2-\alpha}$.
	Note that if $\gamma = 0$ (i.e., $\beta = \frac{1-\alpha}{2-\alpha}$), we have $\x^*_{0} = \x^*$.
	Let $\gamma_s = \beta_s - \frac{1-\alpha_s}{2-\alpha_s}$.
	Note that $\gamma_s \in [0, \beta_s - 1/2)$ for all $s$ and $\lim_{s\to \infty} \gamma_s = 0$
    Let $\tilde{A}_{\alpha,\beta}$ be the matrix defined in
	Theorem~\ref{thm:fLSQ},
	and $\kappa_{\alpha,\beta}$ be its condition number.

	It follows from Theorem~\ref{thm:fLSQ}
	that for $s = 1,\dots$, 
	\begin{equation*}
	    \|\x_s^{(k_s)} - \x^*_{\gamma_s}\|^2 \le r_s^{k_s} \|\x_s^{(0)} - \x^*_{\gamma_s}\|^2
	    \quad
	    \text{where}
	    \quad 
	    r_s = 1- \frac{\eta}{\kappa_{\alpha_s,\beta_s}}.
	\end{equation*}
	Let $\mathcal{E}_s = \|\x_{s}^{(0)} - \x^*_{\gamma_{s}}\|$,
	$e_s = \|\x^*_{\gamma_{s}} - \x^*_{\gamma_{s+1}}\|$
	and $R_s = r_s^{k_s/2}$.
	Since $\x_{s-1}^{(k_{s-1})} = \x_{s}^{(0)}$ for all $s$, we have 
	$\mathcal{E}_s \le R_{s-1}\mathcal{E}_{s-1} + e_{s-1}$.
	Also, observe that 
	$\|\x_s^{(k_s)} - \x^*_{0}\|
	    \le R_s \mathcal{E}_s + \|\x^*_{\gamma_s} - \x^*_{0}\|$.
	By recursively applying $\mathcal{E}_s \le R_{s-1}\mathcal{E}_{s-1} + e_{s-1}$,
	we have 
	$R_{s}\mathcal{E}_s \le 
	    \sum_{k=1}^{s} \left(\prod_{j=0}^{k-1} R_{s-j}\right) e_{s-k}$,
	where $e_0 = \mathcal{E}_1$.
	This gives 
	$\|\x_s^{(k_s)} - \x^*_{0}\|
	    \le \sum_{k=1}^{s} \left(\prod_{j=0}^{k-1} R_{s-j}\right) e_{s-k} + \|\x^*_{\gamma_s} - \x^*_{0}\|$.
	Observe that for any $\gamma, \gamma' \in [0, \beta-1/2)$,
    \begin{align*}
    &\x^*_{\gamma'} - \x^*_{\gamma}
    \\
    &= \left[\left(WW^\top + \gamma' RR^\top\right)^{-1} - \left(WW^\top + \gamma RR^\top\right)^{-1}\right] W(y-W^\top \bc)
    \\
    &= (\gamma - \gamma')\left(WW^\top + \gamma' RR^\top\right)^{-1}
    RR^\top 
    \left(WW^\top + \gamma RR^\top\right)^{-1}
    WW^\top(\x^*-\bc).
    \end{align*}
    Assuming $B_{\max} = \sup_{\gamma \in [0,\beta-1/2)}\|\left(WW^\top + \gamma RR^\top\right)^{-1}\|$ is finite, we obtain  
    $\|\x^*_{\gamma'} - \x^*_{\gamma}\| \le
    C|\gamma - \gamma'|$,
    where $C = B_{\max}^2\|W\|^2\|R\|^2\|\x^* - \bc\|$.
    Therefore, $\|\x^*_{\gamma_s} - \x^*_{0}\| \le C|\gamma_s|$
    and $e_s \le C |\gamma_{s} - \gamma_{s+1}|$,  $\forall s \ge 1$,
    which gives 
    \begin{align*}
	    &\|\x_s^{(k_s)} - \x^*\|
	    \\
	    &\le \left(\prod_{j=0}^{s-1} R_{s-j}\right)
	    \|\x^{(0)} - \x^*_{\gamma_{1}}\|
	    + C\left\{ \sum_{k=1}^{s-1} \left(\prod_{j=0}^{k-1} R_{s-j}\right) |\gamma_{s-k} - \gamma_{s-k+1}| + |\gamma_s|
	    \right\}.
	\end{align*}
\end{proof}

\section{Complexity Calculations} \label{app:complexity}
Firstly, we note that given $\x$, 
the evaluation of the difference between the network prediction and the output data (misfit)
costs $\mathcal{O}(mn)$ FLOPS as 
\begin{align*}
    \text{MISFIT}_{\x}:= W^\top \bm{a}_3 - \bm{y},
\end{align*}
where 
$(W)_{ji} = \phi(a_{1,j}z_i + a_{2,j})$,
$(\bm{y})_i = y_i$,
and $(\bm{a}_3)_j = a_{3, j}$.
For $j=1,\dots, n$, let
\begin{align*}
    \text{MISFIT}_{t,\x}^{[i]}(u) &= 
    \sum_{l=1,l\ne j}^n a_{3,l}\phi(a_{1,l}z_i+a_{2,l}) - y_i + a_{3,j}\phi(uz_i + a_{2,j}), \text{ if } t = j, \\
    \text{MISFIT}_{t,\x}^{[i]}(u) &= \sum_{l=1,l\ne j}^n a_{3,l}\phi(a_{1,l}z_i+a_{2,l}) - y_i + a_{3,j}\phi(a_{1,j}z_i + u),
    \text{ if } t = n + j, \\
    \text{MISFIT}_{t,\x}^{[i]}(u) &= \sum_{l=1,l\ne j}^n a_{3,l}\phi(a_{1,l}z_i+a_{2,l}) - y_i + u\phi(a_{1,j}z_i + a_{2,j}),
    \text{ if } t = 2n + j.
\end{align*}
Observe that for $j=1,\dots, n$,
\begin{align*}
    \text{MISFIT}_{t,\x}^{[i]}(u)
    =
    \begin{cases}
     (\text{MISFIT}_{\x})_i + a_{3,j}(\phi(uz_i+a_{2,j}) - \phi(a_{1,j}z_i+a_{2,j})) & \text{if } t=j, \\
     (\text{MISFIT}_{\x})_i + a_{3,j}(\phi(a_{1,j}z_i+u) - \phi(a_{1,j}z_i+a_{2,j})) & \text{if } t=n+j, \\
     (\text{MISFIT}_{\x})_i + (u-a_{3,j})\phi(a_{1,j}z_i+a_{2,j}) & \text{if } t=2n+j.
    \end{cases}
\end{align*}
This shows that
if $\text{MISFIT}_{\x}$ and $W$ are already computed and stored,
the function $\text{MISFIT}_{t,\x}^{[i]}(u)$
can be evaluated with almost no computational cost.
Assuming evaluation of $\phi$ takes 1 FLOPS, 
a single evaluation of $\text{MISFIT}_{t,\x}^{[i]}(u)$ costs
at most 6 FLOPS.
With the function $\text{MISFIT}_{t,\x}^{[i]}(u)$,
$f_{j,\x}'$ and $f_{j,\x}''$ are given as follows:
For $t = j$ where $j = 1,\dots, n$, 
\begin{align*}
    f'_{t,\x}(u)
    &= \sum_{i=1}^m
    \text{MISFIT}_{t,\x}^{[i]}(u)
    a_{3,j}\phi'(uz_i+a_{2,j})z_i,
    \\
    f''_{t,\x}(u)
    &= \sum_{i=1}^m
    \left\{\text{MISFIT}_{t,\x}^{[i]}(u)
    a_{3,j}\phi''(uz_i+a_{2,j})z_i^2
    +
    (a_{3,j}\phi'(uz_i+a_{2,j})z_i)^2\right\}.
\end{align*}
For $t =n+j$ where $j=1,\dots,n$,
\begin{align*}
    f'_{t,\x}(u)
    &= \sum_{i=1}^m
    \text{MISFIT}_{t,\x}^{[i]}(u)
    a_{3,j}\phi'(a_{1,j}z_i+u),
    \\
    f''_{t,\x}(u)
    &= \sum_{i=1}^m
    \left\{\text{MISFIT}_{t,\x}^{[i]}(u)
    a_{3,j}\phi''(a_{1,j}z_i+u)
    +
    (a_{3,j}\phi'(a_{1,j}z_i+u))^2\right\}.
\end{align*}
For $t =2n+j$ where $j=1,\dots,n$,
\begin{align*}
    f'_{t,\x}(u)
    &= \sum_{i=1}^m
    \text{MISFIT}_{t,\x}^{[i]}(u)\phi(a_{1,j}z_i+a_{2,j}),
    \quad
    f''_{t,\x}(u)
    = \sum_{i=1}^m
    (\phi(a_{1,j}z_i+a_{2,j}))^2.
\end{align*}
Assuming the evaluations of $\phi, \phi', \phi''$ take 1 FLOPS each,
it can be checked that 
a single evaluation of $f_{t,\x}'(u)$/$f_{t,\x}''(u)$
takes at most 
$13m$/$20m$ FLOPS.
	
	\bibliographystyle{siamplain}
	\bibliography{references}
\end{document}